\documentclass[12pt]{article}
\usepackage[utf8]{inputenc}
\usepackage{float}

\usepackage{soul} %
\usepackage{setspace} 
\usepackage[most]{tcolorbox}
\usepackage{shellesc}
\usepackage{mathrsfs,amsthm,xcolor,verbatim,bbm,amsmath,amsfonts,amssymb,nicefrac,enumitem,hyperref,bm,mathtools,xparse,etoolbox,textcomp}
\usepackage{cite}
\usepackage[outer=2cm, inner=3cm,top=4cm,bottom=5cm, twoside, a4paper]{geometry}
\usepackage[capitalise,nameinlink]{cleveref} 
\usepackage{xurl}
\usepackage{frontespizio}
\usepackage{caption, booktabs} %
\usepackage{breakurl}

\usepackage{tikz}
\usetikzlibrary{matrix, arrows.meta}
\usepackage{tikz-3dplot}
\usepackage{pst-solides3d}
\usetikzlibrary{automata,chains}
\usepackage{subcaption}

\tcbuselibrary{breakable,skins}

\usepackage[textwidth=0.7in]{todonotes}

\crefname{enumi}{item}{items}

\crefname{equation}{}{}
\crefname{subsection}{Subsection}{Subsections}

\definecolor{darkblue}{rgb}{0,0,0.75}
\definecolor{darkgreen}{rgb}{0,0.75,0}
\definecolor{darkred}{rgb}{0.75,0,0}
\hypersetup{colorlinks=true,linkcolor=darkblue,citecolor=darkgreen,urlcolor=darkblue}

\theoremstyle{plain}
\newtheorem{theorem}{Theorem}[section]

\newtheorem{prop}[theorem]{Proposition}

\newtheorem{setting}[theorem]{Setting}

\newtheorem{definition}[theorem]{Definition}

\theoremstyle{remark}
\newtheorem{remark}[theorem]{Remark}

\theoremstyle{definition}

\DeclareFontEncoding{LS1}{}{}
\DeclareFontSubstitution{LS1}{stix}{m}{n}
\DeclareMathAlphabet{\mathscr}{LS1}{stixscr}{m}{n}

\newcommand{\E}{\mathbb{E}}
\newcommand{\Pp}{\mathbb{P}}

\newcommand{\R}{\mathbb{R}}
\newcommand{\N}{\mathbb{N}}

\renewcommand{\c}[1]{\mathfrak{c}^{#1}}

\DeclarePairedDelimiterX{\infdivx}[2]{(}{)}{%
  #1\;\delimsize\|\;#2%
}

\newcommand{\m}{\mathbf{m}}

\newcommand{\smallsum}{\textstyle\sum}

\newcommand{\Inn}{U}

\newcommand{\cgamma}{\Gamma}
\newcommand{\g}{g} 
\newcommand{\f}{g} 
\newcommand{\gdet}{g}

\newcommand{\Gbatch}{G}

\newcommand{\const}{\mathfrak{c}}

\newcommand{\dimm}{\delta}

\newcommand{\C}{\mathfrak{C}}

\newcommand{\cF}{\mathcal{F}}

\newcommand{\fd}{\mathfrak{d}}

\makeatletter
\newcommand{\@normibar}{\vrule\@width 1.2\p@}
\newcommand{\normi}{\@ifstar\@xnormi\@normi}
\newcommand{\@xnormi}[1]{%
  \left.\kern-\nulldelimiterspace
  \@normibar
  #1
  \@normibar
  \right.\kern-\nulldelimiterspace
}
\newcommand\@normi[2][]{%
  \mathopen{\vphantom{#1|}\mkern2mu\@normibar\mkern2mu}
  #2
  \mathclose{\vphantom{#1|}\mkern2mu\@normibar\mkern2mu}
}
\makeatother

\newcommand{\qandq}{\qquad\text{and}\qquad}
\newcommand{\andq}{\text{and}\qquad}

\usepackage{cleveref}

\ExplSyntaxOn

\seq_new:N \g_cflist_loaded
\seq_new:N \g_cflist_pending

\NewDocumentCommand{\cfadd} { m } {
  \seq_if_in:NnF \g_cflist_loaded { #1 } {
    \seq_if_in:NnF \g_cflist_pending { #1 } {
      \seq_gput_right:Nn \g_cflist_pending { #1 }
    }
  }
}

\NewDocumentCommand{\cfconsiderloaded} { m } {
  \seq_gput_right:Nn \g_cflist_loaded {#1}
}

\NewDocumentCommand{\cfremove} { m } {
  \seq_gremove_all:Nn \g_cflist_pending { #1 }
}

\NewDocumentCommand{\cfload} { o } {
  \seq_if_empty:NTF \g_cflist_pending {
    \IfValueTF{#1}{\ignorespaces}{\unskip}
  } {
    (cf.\ \cref{\seq_use:Nn \g_cflist_pending {,}})\IfValueTF{#1}{#1~}{\unskip}
    \seq_gconcat:NNN \g_cflist_loaded \g_cflist_loaded \g_cflist_pending
    \seq_gclear:N \g_cflist_pending
    \IfValueT{#1}{\ignorespaces}
  }
}

\NewDocumentCommand{\cfclear} {} {
  \seq_gclear:N \g_cflist_loaded
  \seq_gclear:N \g_cflist_pending
}

\NewDocumentCommand{\cfout} { o } {
  \seq_if_empty:NTF \g_cflist_pending {\unskip\IfValueT{#1}{\ignorespaces}} {
    (cf.\ \cref{\seq_use:Nn \g_cflist_pending {,}})\IfValueTF{#1}{#1~}{\unskip}
    \seq_gclear:N \g_cflist_pending
    \IfValueT{#1}{\ignorespaces}
  }
}

\NewDocumentCommand{\ifnocf} { m } {
  \seq_if_empty:NT \g_cflist_pending { #1 }
}

\ExplSyntaxOff

\ExplSyntaxOn

\bool_new:N \g_noteobserve

\NewDocumentCommand{\setnote}{}{
  \bool_gset_true:N \g_noteobserve
}

\NewDocumentCommand{\setobserve}{}{
  \bool_gset_false:N \g_noteobserve
}

\NewDocumentCommand{\nobs}{ o }{
  \IfValueT{#1}{
    \str_if_eq:noTF {note} {#1} {
      \bool_gset_true:N \g_noteobserve
    } {
      \str_if_eq:noTF {Note} {#1} {
        \bool_gset_true:N \g_noteobserve
      } {
        \bool_gset_false:N \g_noteobserve
      }
    }
  }
  \bool_if:nTF { \g_noteobserve } {
    \bool_gset_false:N \g_noteobserve
    note
  } {
    \bool_gset_true:N \g_noteobserve
    observe
  }
  \IfValueF{#1}{~}
}

\NewDocumentCommand{\Nobs}{ o }{
  \IfValueT{#1}{
    \str_if_eq:noTF {note} {#1} {
      \bool_gset_true:N \g_noteobserve
    } {
      \str_if_eq:noTF {Note} {#1} {
        \bool_gset_true:N \g_noteobserve
      } {
        \bool_gset_false:N \g_noteobserve
      }
    }
  }
  \bool_if:nTF { \g_noteobserve } {
    \bool_gset_false:N \g_noteobserve
    Note
  } {
    \bool_gset_true:N \g_noteobserve
    Observe
  }
  \IfValueF{#1}{~}
}

\ExplSyntaxOff

\ExplSyntaxOn

\bool_new:N \g_hencetherefore

\NewDocumentCommand{\hence}{ o }{
  \IfValueT{#1}{
    \str_if_eq:noTF {hence} {#1} {
      \bool_gset_true:N \g_hencetherefore
    } {
      \str_if_eq:noTF {Hence} {#1} {
        \bool_gset_true:N \g_hencetherefore
      } {
        \bool_gset_false:N \g_hencetherefore
      }
    }
  }
  \bool_if:nTF { \g_hencetherefore } {
    \bool_gset_false:N \g_hencetherefore
    hence
  } {
    \bool_gset_true:N \g_hencetherefore
    therefore
  }
  \IfValueF{#1}{~}
}

\NewDocumentCommand{\Hence}{ o }{
  \IfValueT{#1}{
    \str_if_eq:noTF {hence} {#1} {
      \bool_gset_true:N \g_hencetherefore
    } {
      \str_if_eq:noTF {Hence} {#1} {
        \bool_gset_true:N \g_hencetherefore
      } {
        \bool_gset_false:N \g_hencetherefore
      }
    }
  }
  \bool_if:nTF { \g_hencetherefore } {
    \bool_gset_false:N \g_hencetherefore
    Hence,~we~obtain
  } {
    \bool_gset_true:N \g_hencetherefore
    Therefore,~we~obtain
  }
  \IfValueF{#1}{~}
}

\ExplSyntaxOff

\ExplSyntaxOn

\seq_const_from_clist:Nn \g_prove_mru {
  establish,
  demonstrate,
  prove,
  show,
  imply,
  ensure
}

\prop_new:N \l__verbs
\prop_put:Nnn \l__verbs {show} {shows}
\prop_put:Nnn \l__verbs {imply} {implies}
\prop_put:Nnn \l__verbs {demonstrate} {demonstrates}
\prop_put:Nnn \l__verbs {prove} {proves}
\prop_put:Nnn \l__verbs {establish} {establishes}
\prop_put:Nnn \l__verbs {ensure} {ensures}
\prop_put:Nnn \l__verbs {assure} {assures}

\tl_new:N \g_wordtmp
\seq_new:N \l_mytmps

\cs_generate_variant:Nn \str_if_in:nnTF { nVTF }
\cs_generate_variant:Nn \str_if_in:nnTF { xVTF }

\NewDocumentCommand{\prove}{ o }{
  \IfValueTF{#1}{
    \seq_clear:N \l_mytmps
    \seq_map_inline:Nn \g_prove_mru {
      \str_if_eq:nnTF {##1} {ensure} {
        \str_set:Nn \l_temps {n}
      } {
        \str_set:Nx \l_temps {\str_head_ignore_spaces:n {##1}}
      }
      \str_if_in:xVTF {#1} \l_temps {
        \seq_put_right:Nn \l_mytmps {##1}
      } { }
    }
    \seq_get_right:NN \l_mytmps \g_wordtmp
  } {
    \seq_get_right:NN \g_prove_mru \g_wordtmp
  }
  \tl_use:N \g_wordtmp
  \IfValueTF{#1}{}{~}
  \seq_gput_left:NV \g_prove_mru \g_wordtmp
  \seq_gremove_duplicates:N \g_prove_mru
}

\NewDocumentCommand{\proves}{ o }{
  \IfValueTF{#1}{
    \seq_clear:N \l_mytmps
    \seq_map_inline:Nn \g_prove_mru {
      \str_if_eq:nnTF {##1} {ensure} {
        \str_set:Nn \l_temps {n}
      } {
        \str_set:Nx \l_temps {\str_head_ignore_spaces:n {##1}}
      }
      \str_if_in:xVTF {#1} \l_temps {
        \seq_put_right:Nn \l_mytmps {##1}
      } { }
    }
    \seq_get_right:NN \l_mytmps \g_wordtmp
  } {
    \seq_get_right:NN \g_prove_mru \g_wordtmp
  }
  \str_set:NV \l_tmpa_str \g_wordtmp
  \prop_get:NVN \l__verbs \l_tmpa_str \l_tmpa_tl
  \tl_use:N \l_tmpa_tl
  \IfValueTF{#1}{}{~}
  \seq_gput_left:NV \g_prove_mru \g_wordtmp
  \seq_gremove_duplicates:N \g_prove_mru
}

\newcommand{\llabel}[1]{\savelabel{#1}\label{\loc.#1}\ignorespaces}

\clist_new:N \l_localreflist
\clist_new:N \l_reflist

\NewDocumentCommand{\lref} { m } {
  \clist_set:No \l_localreflist {#1}
  \clist_clear:N \l_reflist
  \clist_map_inline:Nn \l_localreflist { \clist_put_right:Nn \l_reflist {\loc.##1} }
  \cref{\l_reflist}
}

\NewDocumentCommand{\Lref} { m } {
  \clist_set:No \l_localreflist {#1}
  \clist_clear:N \l_reflist
  \clist_map_inline:Nn \l_localreflist { \clist_put_right:Nn \l_reflist {\loc.##1} }
  \Cref{\l_reflist}
}

\NewDocumentCommand{\itref}{ m m }{
  \clist_set:No \l_localreflist {#2}
  \clist_clear:N \l_reflist
  \clist_map_inline:Nn \l_localreflist { \clist_put_right:Nn \l_reflist {#1.##1} }
  \cref{\l_reflist}~in~\cref{#1}
}

\seq_new:N \l_enum_seq
\int_new:N \l_num_items

\bool_new:N \g_commaused_bool

\providecommand{\comma}{}

\cs_new:Nn \enum_it:nn {
  \int_case:nnF {\l_num_items - #1} {
    {0} {
      \renewcommand{\comma}{}
      #2\space
    }
    {1} {
      \bool_gset_false:N \g_commaused_bool
      \renewcommand{\comma}{,~\bool_gset_true:N \g_commaused_bool}
      #2
      \bool_if:NTF \g_commaused_bool {} {,~}
      and~
    }
  } {
    \bool_gset_false:N \g_commaused_bool
    \renewcommand{\comma}{,~\bool_gset_true:N \g_commaused_bool}
    #2
    \bool_if:NTF \g_commaused_bool {} {,~}
  }
}

\cs_new:Nn \enum_it_U:nn {
  \int_case:nnF {\l_num_items - #1} {
    {0} {
      \renewcommand{\comma}{}
      #2
      \space
    }
    {1} {
      \bool_gset_false:N \g_commaused_bool
      \renewcommand{\comma}{,~\bool_gset_true:N \g_commaused_bool}
      #2
      \bool_if:NTF \g_commaused_bool {} {,~}
      and~
    }
  } {
    \bool_gset_false:N \g_commaused_bool
    \renewcommand{\comma}{,~\bool_gset_true:N \g_commaused_bool}
    \int_compare:nTF {#1=1} {
      \text_titlecase_first:n {#2}
    } {
      #2
    }
    \bool_if:NTF \g_commaused_bool {} {,~}
  }
}

\cs_new:Nn \uc:n {\MFUsentencecase{#1}}
\cs_generate_variant:Nn \uc:n {x}

\cs_generate_variant:Nn \str_uppercase:n {x}
\cs_generate_variant:Nn \tl_if_eq:nnTF {onTF}

\cs_new:Nn \enum:nnnn {
  \seq_set_split:Nnn \l_enum_seq ; {#1}
  \seq_remove_all:Nn \l_enum_seq { }
  \seq_remove_all:Nn \l_enum_seq {#2}
  \seq_log:N \l_enum_seq
  \int_set:Nn \l_num_items {\seq_count:N \l_enum_seq}
  \int_log:N \l_num_items
  \int_case:nnF {\l_num_items} {
    { 0 } { 0 }
    { 1 } {
      \IfBooleanTF{#4} {
        \tl_set:Nn \l_text_case_exclude_arg_tl {\cref}
        \text_titlecase_first:n {\seq_use:Nn \l_enum_seq {}}
      } {
        \seq_use:Nn \l_enum_seq {}
      }
      \space
      \tl_if_eq:onTF{#3}{-}{}{
        \bool_if:NTF \l_plural_bool {
          \prove[#3]~
        } {
          \proves[#3]~
        }
      }
    }
    { 2 } {
      \IfBooleanTF{#4} {
        \tl_set:Nn \l_text_case_exclude_arg_tl {\cref}
        \text_titlecase_first:n {\seq_use:Nn \l_enum_seq {~and~}}
      } {
        \seq_use:Nn \l_enum_seq {~and~}
      }
      \space
      \tl_if_eq:onTF{#3}{-}{}{
        \prove[#3]~
      }
    }
  } {
    \IfBooleanTF{#4} {
      \tl_set:Nn \l_text_case_exclude_arg_tl {\cref}
      \seq_indexed_map_function:NN \l_enum_seq \enum_it_U:nn
    } {
      \seq_indexed_map_function:NN \l_enum_seq \enum_it:nn
    }
    \tl_if_eq:onTF{#3}{-}{}{
      \prove[#3]~
    }
  }
}

\cs_generate_variant:Nn \enum:nnnn {nxnn}
\cs_generate_variant:Nn \enum:nnnn {nxxn}

\NewDocumentCommand{\enum}{O{} m O{-} s}{
  \IfBooleanTF{#4}{
    \enum:nxnn {#2} {#1} {sindep} \BooleanFalse
  } {
    \enum:nxxn {#2} {#1} {#3} \BooleanFalse
  }
}

\bool_new:N \g_arg_start_bool
\bool_gset_true:N \g_arg_start_bool

\NewDocumentCommand{\startnewargseq}{}{\bool_gset_true:N \g_arg_start_bool \tl_set:Nn \g_label_tl {}}

\cs_generate_variant:Nn \seq_if_in:NnTF {NxTF}
\cs_generate_variant:Nn \seq_remove_all:Nn {Nx}

\int_new:N \l_random_int

\cs_generate_variant:Nn \tl_if_head_eq_catcode:nNTF {oNTF}
\cs_generate_variant:Nn \tl_if_head_eq_catcode:nNTF {VNTF}

\cs_generate_variant:Nn \tl_if_head_eq_catcode:nNTF {eNTF}
\cs_generate_variant:Nn \tl_log:n {o}
\cs_generate_variant:Nn \tl_log:n {f}
\cs_generate_variant:Nn \tl_log:n {x}
\cs_generate_variant:Nn \tl_log:n {e}

\cs_generate_variant:Nn \tl_if_in:nnTF {onTF}
\cs_generate_variant:Nn \tl_if_in:NnTF {NeTF}

\cs_generate_variant:Nn \tl_if_head_eq_meaning:nNTF {VNTF}

\seq_const_from_clist:Nn \g_arg_mru_this {
  Ahpr,
  Tapr,
  Ctapr,
  H
}

\seq_const_from_clist:Nn \g_arg_mru_nothis {
  Ia,
  Nwc,
  N,
  Itns,
  Fm,
  Itnswc,
  Mo
}

\seq_new:N \l_arg_seq
\tl_new:N \l_cons_tl
\tl_new:N \l_dummy_tl

\bool_new:N \g_debug_bool
\bool_gset_false:N \g_debug_bool

\bool_new:N \l_insidearg_bool

\bool_new:N \g_firstargletter_bool

\sys_gset_rand_seed:n {0903}

\bool_new:N \l_plural_bool
\tl_new:N \l_arg_verbs_tl

\NewDocumentCommand{\argument}{mom}{
  \bool_set_false:N \l_plural_bool
  \tl_set:Nn \l_arg_verbs_tl {sindep}
  \keys_define:nn { benno/argument } {
    plural .value_forbidden:n = true,
    plural .code:n = {\bool_set_true:N \l_plural_bool},
    verbs .value_required:n = false,
    verbs .tl_set:N = \l_arg_verbs_tl,
  }
  \IfValueT{#2}{
    \keys_set:nn { benno/argument } {#2}
  }
  \bool_log:N \l_plural_bool
  \bool_gset_true:N \l_insidearg_bool
  \seq_set_split:Nnn \l_arg_seq ; {#1}
  \seq_remove_all:Nn \l_arg_seq { }
  \seq_log:N \l_arg_seq
  \tl_set:Nn \l_cons_tl {#3}
  \tl_trim_spaces:N \l_cons_tl
  \seq_if_in:NxTF \l_arg_seq {\lref{\g_label_tl}} {
    \seq_remove_all:Nx \l_arg_seq {\lref{\g_label_tl}}
    \seq_get_left:NNTF \l_arg_seq \l_dummy_tl {
      \tl_trim_spaces:N \l_dummy_tl
      \bool_gset_false:N \g_firstargletter_bool
      \tl_if_head_eq_catcode:VNTF \l_dummy_tl a {
        \bool_gset_true:N \g_firstargletter_bool
      } {
        \tl_if_head_eq_meaning:VNTF \l_dummy_tl {\cref} {
          \tl_set:Nx \l_tmpa_tl {\tl_tail:N \l_dummy_tl}
          \tl_set:Nx \l_tmpb_tl {\tl_head:N \l_tmpa_tl}
          \bool_gset_true:N \g_firstargletter_bool
          \tl_if_in:NeTF \l_tmpb_tl {lem\c_colon_str} {} {
            \tl_if_in:NeTF \l_tmpb_tl {thm\c_colon_str} {} {
              \tl_if_in:NeTF \l_tmpb_tl {prop\c_colon_str} {} {
                \tl_if_in:NeTF \l_tmpb_tl {cor\c_colon_str} {} {
                  \bool_gset_false:N \g_firstargletter_bool
                }
              }
            }
          }
        } {
        }
      }
      \bool_if:NTF \g_firstargletter_bool {
        \seq_set_eq:NN \l_tmpa_seq \g_arg_mru_this
        \seq_remove_all:Nn \l_tmpa_seq {H}
        \seq_get_right:NN \l_tmpa_seq \l_tmpa_tl
        \int_case:nnF {\seq_count:N \l_arg_seq} {
          {1} {
            \str_case:VnF {\l_tmpa_tl} {
              {Ahpr} {
                \bool_if:NT \g_debug_bool {C1.1}
                \seq_gput_left:Nn \g_arg_mru_this {Ahpr}
                \seq_gremove_duplicates:N \g_arg_mru_this
                \enum:nxnn {#1} {\lref{\g_label_tl}} {-} {\BooleanTrue}
                \hence~
                \bool_if:NTF \l_plural_bool {
                  \prove[\l_arg_verbs_tl]~\ignorespaces #3
                } {
                  \proves[\l_arg_verbs_tl]~\ignorespaces #3
                }
              }
              {Tapr} {
                \bool_if:NT \g_debug_bool {C1.2}
                \seq_gput_left:Nn \g_arg_mru_this {Tapr}
                \seq_gremove_duplicates:N \g_arg_mru_this
                \enum[\lref{\g_label_tl}]{
                  This;
                  #1
                }[\l_arg_verbs_tl]\ignorespaces #3
              }
              {Ctapr} {
                \bool_if:NT \g_debug_bool {C1.3}
                \seq_gput_left:Nn \g_arg_mru_this {Ctapr}
                \seq_gremove_duplicates:N \g_arg_mru_this
                Combining~
                \enum[\lref{\g_label_tl}]{
                  this;
                  #1
                } \proves[\l_arg_verbs_tl]~\ignorespaces #3
              }
            } {}
          }
        } {
          \str_case:VnF {\l_tmpa_tl} {
             {Ahpr} {
              \bool_if:NT \g_debug_bool {C2.1}
              \seq_gput_left:Nn \g_arg_mru_this {Ahpr}
              \seq_gremove_duplicates:N \g_arg_mru_this
              \enum:nxnn {#1} {\lref{\g_label_tl}} {-} {\BooleanTrue}
              \hence~
              \prove[\l_arg_verbs_tl]~\ignorespaces #3
            }
            {Tapr} {
              \bool_if:NT \g_debug_bool {C2.2}
              \seq_gput_left:Nn \g_arg_mru_this {Tapr}
              \seq_gremove_duplicates:N \g_arg_mru_this
              \enum[\lref{\g_label_tl}]{
                This;
                #1
              }[\l_arg_verbs_tl]\ignorespaces #3
            }
            {Ctapr} {
              \int_case:nn {\int_rand:nn {0} {1}} {
                {0} {
                  \bool_if:NT \g_debug_bool {C2.3}
                  \seq_gput_left:Nn \g_arg_mru_this {Ctapr}
                  \seq_gremove_duplicates:N \g_arg_mru_this
                  Combining~
                  \enum[\lref{\g_label_tl}]{
                    this;
                    #1
                  } \proves[\l_arg_verbs_tl]~\ignorespaces #3
                }
                {1} {
                  \bool_if:NT \g_debug_bool {C2.4}
                  \seq_gput_left:Nn \g_arg_mru_this {Ctapr}
                  \seq_gremove_duplicates:N \g_arg_mru_this
                  Combining~
                  \enum:nxnn {#1} {\lref{\g_label_tl}} {-} {\BooleanFalse}
                  \hence~
                  \proves[\l_arg_verbs_tl]~\ignorespaces #3
                }
              }
            }
          } {}
        }
      } {
        \seq_set_eq:NN \l_tmpa_seq \g_arg_mru_this
        \seq_remove_all:Nn \l_tmpa_seq {H}
        \seq_remove_all:Nn \l_tmpa_seq {Ahpr}
        \seq_get_right:NN \l_tmpa_seq \l_tmpa_tl
        \int_case:nnF {\seq_count:N \l_arg_seq} {
          {1} {
            \str_case:VnF {\l_tmpa_tl} {
              {Tapr} {
                \bool_if:NT \g_debug_bool {C3.1}
                \seq_gput_left:Nn \g_arg_mru_this {Tapr}
                \seq_gremove_duplicates:N \g_arg_mru_this
                \enum[\lref{\g_label_tl}]{
                  This;
                  #1
                }[\l_arg_verbs_tl]\ignorespaces #3
              }
              {Ctapr} {
                \bool_if:NT \g_debug_bool {C3.2}
                \seq_gput_left:Nn \g_arg_mru_this {Ctapr}
                \seq_gremove_duplicates:N \g_arg_mru_this
                Combining~
                \enum[\lref{\g_label_tl}]{
                  this;
                  #1
                } \proves[\l_arg_verbs_tl]~\ignorespaces #3
              }
            } {}
          }
        } {
          \str_case:VnF {\l_tmpa_tl} {
            {Tapr} {
              \bool_if:NT \g_debug_bool {C4.1}
              \seq_gput_left:Nn \g_arg_mru_this {Tapr}
              \seq_gremove_duplicates:N \g_arg_mru_this
              \enum[\lref{\g_label_tl}]{
                This;
                #1
              }[\l_arg_verbs_tl]\ignorespaces #3		
            }
            {Ctapr} {
              \int_case:nn {\int_rand:nn {0} {1}} {
                {0} {
                  \bool_if:NT \g_debug_bool {C4.2}
                  \seq_gput_left:Nn \g_arg_mru_this {Ctapr}
                  \seq_gremove_duplicates:N \g_arg_mru_this
                  Combining~
                  \enum[\lref{\g_label_tl}]{
                    this;
                    #1
                  } \proves[\l_arg_verbs_tl]~\ignorespaces #3		
                }
                {1} {
                  \bool_if:NT \g_debug_bool {C4.3}
                  \seq_gput_left:Nn \g_arg_mru_this {Ctapr}
                  \seq_gremove_duplicates:N \g_arg_mru_this
                  Combining~
                  \enum:nxnn {#1} {\lref{\g_label_tl}} {-} {\BooleanFalse}
                  \hence~
                  \proves[\l_arg_verbs_tl]~\ignorespaces #3    
                }
              }
            }
          } {}
        }
      }
    } {
      \tl_if_head_eq_catcode:oNTF \l_cons_tl a {
        \seq_set_eq:NN \l_tmpa_seq \g_arg_mru_this
        \seq_remove_all:Nn \l_tmpa_seq {Ctapr}
        \seq_remove_all:Nn \l_tmpa_seq {Ahpr}
        \seq_get_right:NN \l_tmpa_seq \l_tmpa_tl
        \str_case:VnF {\l_tmpa_tl} {
          {H} {
            \bool_if:NT \g_debug_bool {C5.1}
            \seq_gput_left:Nn \g_arg_mru_this {H}
            \seq_gremove_duplicates:N \g_arg_mru_this
            Hence,~we~obtain~\ignorespaces #3
          }
          {Tapr} {
            \bool_if:NT \g_debug_bool {C5.2}
            \seq_gput_left:Nn \g_arg_mru_this {Tapr}
            \seq_gremove_duplicates:N \g_arg_mru_this
            This~\proves[\l_arg_verbs_tl]~\ignorespaces #3
          }
        } {}
      } {
        \bool_if:NT \g_debug_bool {C6.1}
        \seq_gput_left:Nn \g_arg_mru_this {Tapr}
        \seq_gremove_duplicates:N \g_arg_mru_this
        This~\proves[\l_arg_verbs_tl]~\ignorespaces #3
      }
    } 
  } {
    \int_compare:nNnTF {\seq_count:N \l_arg_seq} = {0} {
      \bool_if:NTF \g_arg_start_bool {
        \bool_if:NT \g_debug_bool {C7.1}
        \Nobs\unskip
        #3
      } {
        \bool_if:NT \g_debug_bool {C7.2}
        \Moreover~
        #3
      }
    } {
      \bool_if:NTF \g_arg_start_bool {
        \bool_if:NT \g_debug_bool {C8.1}
        \tl_log:N \l_arg_verbs_tl
        \Nobs~that~
        \enum{
          #1
        }[\l_arg_verbs_tl]\ignorespaces #3
      } {
        \int_compare:nNnTF {\seq_count:N \l_arg_seq} = {1} {
          \seq_set_eq:NN \l_tmpa_seq \g_arg_mru_nothis
          \seq_remove_all:Nn \l_tmpa_seq {Nwc}
          \seq_remove_all:Nn \l_tmpa_seq {Itnswc}
          \seq_get_right:NN \l_tmpa_seq \l_tmpa_tl
        } {
          \seq_get_right:NN \g_arg_mru_nothis \l_tmpa_tl
        }
        \str_case:VnF {\l_tmpa_tl} {
          {Mo} {
            \bool_if:NT \g_debug_bool {C9.1}
            \seq_gput_left:Nn \g_arg_mru_nothis {Mo}
            \seq_gremove_duplicates:N \g_arg_mru_nothis
            Moreover,~\nobs~that~
            \enum{
              #1
            }[\l_arg_verbs_tl]\ignorespaces #3		
          }
          {Fm} {
            \bool_if:NT \g_debug_bool {C9.2}
            \seq_gput_left:Nn \g_arg_mru_nothis {Fm}
            \seq_gremove_duplicates:N \g_arg_mru_nothis
            Furthermore,~\nobs~that~
            \enum{
              #1
            }[\l_arg_verbs_tl]\ignorespaces #3		
          }
          {Ia} {
            \bool_if:NT \g_debug_bool {C9.3}
            \seq_gput_left:Nn \g_arg_mru_nothis {Ia}
            \seq_gremove_duplicates:N \g_arg_mru_nothis
            In~addition,~\nobs~that~
            \enum{
              #1
            }[\l_arg_verbs_tl]\ignorespaces #3		
          }
          {N} {
            \bool_if:NT \g_debug_bool {C9.4}
            \seq_gput_left:Nn \g_arg_mru_nothis {N}
            \seq_gremove_duplicates:N \g_arg_mru_nothis
            Next,~\nobs~that~
            \enum{
              #1
            }[\l_arg_verbs_tl]\ignorespaces #3		
          }
          {Itns} {
            \bool_if:NT \g_debug_bool {C9.5}
            \seq_gput_left:Nn \g_arg_mru_nothis {Itns}
            \seq_gremove_duplicates:N \g_arg_mru_nothis
            In~the~next~step~we~\nobs~that~
            \enum{
              #1
            }[\l_arg_verbs_tl]\ignorespaces #3		
          }
          {Nwc} {
            \bool_if:NT \g_debug_bool {C9.6}
            \seq_gput_left:Nn \g_arg_mru_nothis {Nwc}
            \seq_gremove_duplicates:N \g_arg_mru_nothis
            Next~we~combine~
            \enum{
              #1
            }to~obtain~\ignorespaces #3
          }
          {Itnswc} {
            \bool_if:NT \g_debug_bool {C9.7}
            \seq_gput_left:Nn \g_arg_mru_nothis {Itnswc}
            \seq_gremove_duplicates:N \g_arg_mru_nothis
            In~the~next~step~we~combine~
            \enum{
              #1
            }to~obtain~\ignorespaces #3
          }
        } {}
      }
    }
  }
  \bool_gset_false:N \g_arg_start_bool
  \bool_gset_false:N \l_insidearg_bool
  \cfload[.]%
  \color{black}
}

\tl_new:N \g_label_tl
\tl_gset:Nn \g_label_tl { }

\NewDocumentCommand{\savelabel}{m}{
  \bool_if:NTF \l_insidearg_bool {
    \tl_gset:Nn \g_label_tl {#1}
  } {
    \tl_gset:Nn \g_label_tl { }
  }
}

\ExplSyntaxOff

\ExplSyntaxOn

\NewDocumentEnvironment {athm} {m m o} {
\str_if_eq:noTF {example} {#1} {
  \bool_gset_true:N \g_example_bool
} {
  \bool_gset_false:N \g_example_bool
}
\cfclear
\begin{samepage}
\begin{tcolorbox}[colback=white!95!gray,
                  colframe=black,
                  boxrule=0.5pt,
                  sharp~corners,
                  enhanced,
                  breakable,
                 ]
\IfNoValueTF{#3}{
\begin{#1}\label{#2}\global\def\loc{#2}
}{
\begin{#1}[#3]\label{#2}\global\def\loc{#2}
}
}{
\end{#1}
\end{tcolorbox}
\end{samepage}
}

\NewDocumentEnvironment {adef} {m} {
\begin{samepage}
\begin{tcolorbox}[colback=white!95!gray,
                  colframe=black,
                  boxrule=0.5pt,
                  sharp~corners,
                  enhanced,
                  breakable,
                 ]
\begin{definition}\label{#1}\global\def\loc{#1}
}{
\end{definition}
\end{tcolorbox}
\end{samepage}
}

\NewDocumentEnvironment{aproof} {} {
\bool_if:NTF \g_example_bool {
  \bool_gset_true:N \g_arg_start_bool
  \begin{proof}[Proof~for~\cref{\loc}]
} {
  \bool_gset_true:N \g_arg_start_bool
  \begin{proof}[Proof~of~\cref{\loc}]
}
\bool_gset_false:N \g_finishproof_bool
}{
\bool_if:NTF \g_finishproof_bool {}
{\finishproofthus}
\end{proof}
}

\NewDocumentCommand{\finishproofthus} {} {
  \bool_gset_true:N \g_finishproof_bool 
  \bool_if:NTF \g_example_bool {
    The~proof~for~\cref{\loc}~is~thus~complete.
  } {
    The~proof~of~\cref{\loc}~is~thus~complete.
  }
}
\NewDocumentCommand{\finishproofthis} {} {
  \bool_gset_true:N \g_finishproof_bool 
  \bool_if:NTF \g_example_bool {
    This~completes~the~proof~for~\cref{\loc}.
  } {
    This~completes~the~proof~of~\cref{\loc}.
  }
}

\ExplSyntaxOff

\NewDocumentEnvironment{cproof}{m}
{\begin{proof}[Proof of \cref{#1}]}%
{\noindent The proof of \cref{#1} is thus complete.
\end{proof}}

\NewDocumentEnvironment{cproof2}{m}
{\begin{proof}[Proof of \cref{#1}]}%
{\noindent This completes the proof of \cref{#1}.
\end{proof}}

\ExplSyntaxOn

\int_new:N \g_furthermore

\NewDocumentCommand{\Moreover}{ o o }{
  \IfValueT{#1}{
    \str_case:nn {#1} {
      {Next} {\int_gset:Nn {\g_furthermore} {0}}      
      {Furthermore} {\int_gset:Nn {\g_furthermore} {1}}
      {Moreover} {\int_gset:Nn {\g_furthermore} {2}}
      {In~addition} {\int_gset:Nn {\g_furthermore} {3}}
      {note} {\bool_gset_true:N \g_noteobserve}
      {observe} {\bool_gset_false:N \g_noteobserve}
    }
    \IfValueT{#2}{
      \str_case:nn {#2} {
        {Next} {\int_gset:Nn {\g_furthermore} {0}}        
        {Furthermore} {\int_gset:Nn {\g_furthermore} {1}}
        {Moreover} {\int_gset:Nn {\g_furthermore} {2}}
        {In~addition} {\int_gset:Nn {\g_furthermore} {3}}
        {note} {\bool_gset_true:N \g_noteobserve}
        {observe} {\bool_gset_false:N \g_noteobserve}
      }
    }
  }
  \int_case:nn { \int_mod:nn {\g_furthermore} {4} } {
    { 0 } { Next~\nobs that}    
    { 1 } { Furthermore,~\nobs that}
    { 2 } { Moreover,~\nobs that}
    { 3 } { In~addition,~\nobs that}
  }
  \int_incr:N \g_furthermore
  \peek_charcode:NTF , {  } { ~ }
}

\ExplSyntaxOff

\hypersetup{
    colorlinks,
    linkcolor={blue!80!black},
    citecolor={green},
    urlcolor={blue!80!black}
}

\makeatletter
\ExplSyntaxOn
\seq_new:N \g_abbrs
\prop_new:N \g_abbr_counts
\tl_new:N \l_abbr_count_tl

\ExplSyntaxOn

\bool_new:N \g_forexample

\NewDocumentCommand{\eg}{ o }{
	\IfValueT{#1}{
		\str_if_eq:noTF {fe} {#1} {
			\bool_gset_true:N \g_forexample
		} {\bool_gset_false:N \g_forexample}
	}
	\bool_if:nTF { \g_forexample } {
		\bool_gset_false:N \g_forexample
		for~example
	}{
		\bool_gset_true:N \g_forexample
		for~instance
	}
}

\NewDocumentCommand{\abbr}{m m O{#1} m m O{#4} m}{
	\expandafter\newcommand\csname#3\endcsname[1][]{
		\seq_if_in:NnTF \g_abbrs {#1} {
			\prop_get:NnN \g_abbr_counts {#1} \l_abbr_count_tl
			\prop_gput:Nnx \g_abbr_counts {#1} {\int_eval:n {\l_abbr_count_tl + 1}}
			\hyperref[#1]{#7}
		} {
			\seq_gput_left:Nn \g_abbrs {#1}
			\prop_gput:Nnn \g_abbr_counts {#1} {1}
			\expandafter\gdef\csname#1@def\endcsname{#2}
			\phantomsection\label{#1}
			\str_if_eq:nnTF{##1}{}{\emph{#2}}{##1}~(\hyperref[#1]{#7})
		}
	}
	\expandafter\newcommand\csname#6\endcsname[1][]{
		\seq_if_in:NnTF \g_abbrs {#1} {
			\prop_get:NnN \g_abbr_counts {#1} \l_abbr_count_tl
			\prop_gput:Nnx \g_abbr_counts {#1} {\int_eval:n {\l_abbr_count_tl + 1}}
			\hyperref[#1]{#4}
		} {
			\expandafter\gdef\csname#1@def\endcsname{#5}
			\seq_gput_left:Nn \g_abbrs {#1}
			\prop_gput:Nnn \g_abbr_counts {#1} {1}
			\phantomsection\label{#1}
			\str_if_eq:nnTF{##1}{}{\emph{#5}}{##1}~(\hyperref[#1]{#4})
		}
	}
}

\ExplSyntaxOff
\makeatother

\ExplSyntaxOff
\makeatother
\abbr{SGD}{Stochastic gradient descent}{SGDs}{stochastic gradient descents}{SGD}
\abbr{iid}{independent and identically distributed}{i.i.d}{independent and identically distributed}{i.i.d.}
\abbr{AI}{artificial intelligence}{AIs}{artificial intelligence}{AI} 
\abbr{SAP}{stochastic approximation problem}{SAPs}{stochastic approximation problem}{SAP} 
\abbr{SOP}{stochastic optimization problem}{SOPs}{stochastic optimization problem}{SOP} 
\abbr{LLM}{large language model}{LLMs}{large language models}{LLM} 
\abbr{Adam}{adaptive moment estimation}{Adams}{adaptive moment estimation}{Adam} 
\abbr{AdamW}{\Adam\ with decoupled weight decay}{AdamWs}{\Adam\ with decoupled weight decay}{AdamW} 
\abbr{MUON}{momentum orthogonalized}{MUONs}{momentum orthogonalizeds}{MUON}

\title{Strong error analysis for the\\ stochastic momentum optimizer}
\author{Davide Gallon$^{1}$ and Arnulf Jentzen$^{2,3}$
\bigskip
    \\
	\small{$^1$Applied Mathematics: Institute for Analysis and Numerics,}
	\\
    \small{University of M\"unster, Germany, e-mail: davide.gallon@uni-muenster.de}
	\smallskip
	\\
	\small{$^2$School of Data Science and School of Artificial Intelligence,}
	\\
	\small{The Chinese University of Hong Kong, Shenzhen (CUHK-Shenzhen),}
    \\
    \small{China, e-mail: ajentzen@cuhk.edu.cn}
	\smallskip
	\\
	\small{$^3$Applied Mathematics: Institute for Analysis and Numerics,}
	\\
	\small{University of M\"unster, Germany, e-mail: ajentzen@uni-muenster.de} 
	\smallskip
	\\
}
\date{\today}

\begin{document}

\maketitle

\begin{abstract}
    \noindent
    Stochastic gradient descent (SGD) optimization schemes are the methods of choice for the optimization of deep neural networks (DNNs) in artificial intelligence (AI) systems. Often not the standard SGD method is used but instead suitable accelerated, adaptive, and/or normalized variants of standard SGD such as Adam, AdamW, and MUON are employed to train large scale AI systems in practically relevant settings. The acceleration (higher order convergence speed) in all these popular optimizers relies on the momentum SGD optimizer. In this work we provide a rigorous error analysis for the momentum SGD optimizer. In particular, we establish convergence rates for the momentum optimizer in terms of the size of the learning rate (step size), the size of the mini-batch, and the size of the one-point convexity constant. 
\end{abstract}

\pagebreak

\tableofcontents

\section{Introduction}

\SGD\ optimization schemes are nowadays the methods of choice to train large scale \AI\ systems \cite{ruder2017overviewgradientdescentoptimization,jentzen2023mathematical}. Instead of the plain vanilla standard \SGD\ method \cite{d8d62392-9a37-31e7-ad3b-37a6f6ee8ef6}, one often employs sophisticated \SGD\ methods involving acceleration, adaptivity, and/or normalization techniques, such as in the popular \Adam\ \cite{kingma2014adam}, the \AdamW\ \cite{loshchilov2017decoupled}, or the \MUON\ by Newton-Schulz \cite{jordan2024muon} optimizers. In particular, according to the \emph{artificial analysis intelligence index} \cite{ArtificialAnalysis}, the smartest \LLMs\ are, where this information has been publicly disclosed, all trained with \AdamW, \MUON, or variants of \MUON\  (\MUON\ Clip, \MUON\ Split). 

The acceleration (higher order convergence speed) for each of these optimizers, in turn, relies on the \emph{momentum \SGD\ optimizer}, which is incorporated into each of these optimization methods. In the situation of deterministic optimization problems it is already known for a long time that the momentum method converges \emph{with a strictly higher rate of convergence} than the standard gradient descent method; see \cite{POLYAK19641} and, \eg, \cite[ Items (i) and (ii) in Theorem 1.2]{dereich2025sharp}. 

\subsection{Main result: Error analysis for momentum stochastic gradient descent (SGD)}

Taking the relevance of the acceleration due to the momentum \SGD\ method into account, in this work we aim to take a closer look at the momentum optimizer and, in particular, we provide an error analysis with convergence rates for the momentum method.
 Specifically, in \cref{theorem1} below we establish an upper bound for the mean square distance $\E\bigl[\| \Theta^J_n - \vartheta \|^2 \bigr]$ between the solution $\vartheta \in \R^\fd$ of the \SAP\ under consideration and the momentum \SGD\ optimization process $(\Theta^J_n )_{ n \in \N_0 }$ with mini-batch size $J \in \N = \{ 1, 2, 3, ... \}$ after $n \in \N_0 = \{0, 1, 2, 3, ...\}$ gradient steps. The natural number $\delta \in \N$ in \cref{theorem1} represents the dimensionality of the data in the considered \SAP\ and the \iid\ random variables $\Inn_{n,m} \colon \Omega \to \R^\dimm $, $(n,m) \in \N_0^2$, in \cref{theorem1} represent the data in the considered \SAP.

\begin{samepage}
\begin{tcolorbox}[colback=white!95!gray,
                  colframe=black,
                  boxrule=0.5pt,
                  sharp corners,
                  enhanced,
                  breakable,
                 ]
\begin{theorem}[\textcolor{red}{Momentum \SGD}]\label{theorem1}
Let $ ( \Omega, \mathcal{F}, \Pp)$ be a probability space, 
let  $\fd, \dimm \in \N$, let $\Inn_{n,m} \colon \Omega \to \R^\dimm $, $(n,m) \in \N_0^2$, be \iid\ random variables,
let $\Gbatch \colon \R^\fd \times \R^\dimm \to \R^\fd$ be measurable, 
let $ \xi, \vartheta \in \R^\fd $, $ \alpha \in [0,1)$, $ c,\cgamma \in (0,\infty)$,  $\lambda \in (0,1)$,
for every $J\in \N$ let 
$ \Theta^J \colon \N_0 \times \Omega \to \R^\fd $ and $ \m^J \colon \N_0 \times \Omega \to \R^\fd $
satisfy for all $ n \in \N $ that 
\begin{equation}\label{theorem1:Theta}
      \textcolor{magenta}{\Theta_0^J = \xi}, \qquad  
      \textcolor{magenta}{\m_0^J=0}, \qquad  \textcolor{magenta}{\Theta_n^J
  = \Theta_{ n - 1 }^J - \cgamma n^{-\lambda} \m_n^J},
\end{equation}
\begin{equation}\label{theorem1:m}
  \andq \textcolor{magenta}{\m_n^J=\alpha \m_{n-1}^J+(1-\alpha)\left[\frac{1}{J} \sum_{j=1}^{J} \Gbatch(\Theta_{n-1}^J, \Inn_{n,j})\right]},
\end{equation}
and assume for all $\theta\in  \R^\fd$ that $\sup_{J,n\in \N }\E \bigl[\|\Gbatch(\theta,\Inn_{0,1})\| +\|\Gbatch(\Theta_{n-1}^J, \Inn_{0,1})\|^2 \bigr] <\infty$ and
\begin{equation}\label{theorem1:coercivity}
    \textcolor{magenta}{\langle\theta-\vartheta, \E[ \Gbatch (\theta, \Inn_{1,1})]\rangle \geq   c \max\{\|\theta-\vartheta\|^2, \|\E[\Gbatch(\theta,\Inn_{1,1})]\|^2\}}.
\end{equation}
Then there exists $\textcolor{magenta}{\C \in (0,\infty)}$ such that for all $\textcolor{magenta}{J,n \in \N}$ it holds that
\begin{equation}
  \textcolor{magenta}{\E\bigl[\|\Theta_n^J -\vartheta \|^2\bigr] \leq 
           \C \bigl(\exp (-\C^{-1} n^{1-\lambda}) + J^{-1} n^{-\lambda} \bigr)}.
    \end{equation}
\end{theorem}
\end{tcolorbox}
\end{samepage}
\cref{theorem1} is an immediate consequence of \cref{cor:main-result}. 
In the following we illustrate the conclusion and some of the assumptions of \cref{theorem1} in words.

The natural number $\fd \in \N$ in \cref{theorem1} represents the dimensionality (the number of degrees of freedom) of the considered \SAP. 
The assumption in \cref{theorem1}  that for all $\theta \in \R^\fd$ it holds that 
$\E [\|\Gbatch(\theta,\Inn_{1,1})\|]< \infty$ (cf.\ \cref{theorem1:m}--\cref{theorem1:coercivity})
ensures that there exists a function $\f \colon \R^{ \fd } \to \R^{ \fd } $
which satisfies for all $\theta \in \R^{ \fd }$ that $\f( \theta ) = \E[ \Gbatch( \theta, \Inn_{ 1, 1 } ) ]$. 
In \cref{theorem1:coercivity} we assume that for all $\theta \in \R^{ \fd }$ it holds that $\langle \theta - \vartheta, \f( \theta )\rangle \geq c \max\{ \| \theta - \vartheta \|^2, \| \f( \theta ) \|^2 \}$ (one-point convexity combined with growth bound). The Cauchy--Schwarz inequality hence ensures that for all $\theta \in \R^\fd$ it holds that
$\langle \theta - \vartheta, \f( \theta ) \rangle \geq c \| \theta - \vartheta \|^2$ (one-point convexity-type assumption)
and 
$\| \f( \theta ) \| \leq c^{ - 1 } \| \theta - \vartheta \|$ (growth bound). 
In the literature these assumptions have been frequently employed for error analyses of 
gradient based optimization methods (cf., \eg, \cite[Corollary 6.1.11 and Lemma 6.1.21]{jentzen2023mathematical}). 

Generally, we think of the function $
\Gbatch = ( \Gbatch(\theta,x) )_{ (\theta,x) \in \R^{\fd}\times\R^\dimm  } \colon \R^\fd \times \R^\dimm \to \R^\fd$ in \cref{theorem1}
as the gradient with respect to the $\theta$-variable of another function $L \allowbreak $ $= \allowbreak( L(\theta,x) )_{  (\theta,x) \in \R^{\fd}\times\R^\dimm } \colon \R^\fd \times \R^\dimm \to [0,\infty)$, 
referred to as \emph{loss function}, and we think of the function $\gdet \colon \R^\fd \to \R^\fd$
as the gradient of another function $\ell \colon \R^\fd  \to [0,\infty)$ satisfying for all $\theta \in \R^{\fd}$ that $\ell( \theta ) = \E\bigl[ L( \theta, \Inn_{ 1, 1 } ) \bigr]$, the objective function of a \SOP, 
but it is not a strict requirement of \cref{theorem1}  that these anti-derivative functions do exist. 

In \cref{theorem1}  the momentum \SGD\ optimization process $\Theta^J \colon \N_0 \times \Omega \to \R^\fd$ is described in \cref{theorem1:Theta}, where the natural number $J \in \N$ denotes the size of the mini-batch. The stochastic process $\m^J \colon \N_0 \times \Omega \to \R^\fd$ represents the auxiliary momentum process appearing in the momentum \SGD\ dynamics. 

\cref{theorem1} bounds the mean square distance $\E[ \| \Theta^J_n - \vartheta \|^2 ]$ 
between the solution $\vartheta \in \R^{ \fd }$ of the \SAP\ under consideration 
and the momentum \SGD\ optimization process $( \Theta^J_n )_{ n \in \N_0 }$ with mini-batch size $J \in \N$ 
after $n \in \N$ gradient steps from above by a possibly arbitrarily large constant $\C$ multiplied by
the sum of the learning rate $\cgamma n^{ -\lambda }$ divided by the size of the mini-batch $J$ 
and the reciprocal of the exponential of a possibly arbitrarily small constant 
multiplied by the training time $\sum_{ k = 1 }^n \cgamma k^{ - \lambda }$.

\subsection{Literature review}\label{lite}
\SGD\ and its variants have been extensively studied in the literature under several algorithmic
parameterizations and convergence criteria.

For plain \SGD, \cite{jentzen2021strong} derives strong \(L^p\) convergence rates under convexity-type assumptions, while \cite{jentzen2020lower} establishes essentially matching mean square lower and upper bounds for polynomially decaying learning rates in a quadratic model. We refer to \cite[Chapter~7]{jentzen2023mathematical} for a broader overview.
For deterministic heavy-ball with constant parameters, \cite{ghadimi2015global} proves an $\mathcal{O}(n^{-1})$ objective-gap bound for the average of all iterates up to step $n$ under smooth convexity and linear convergence under strong convexity. \cite{sun2018nonergodic} establishes instead an $\mathcal{O}(n^{-1})$ bound for the individual iterates under the additional coercivity assumption and linear convergence under restricted strong convexity.

Turning to stochastic momentum methods, we refer, \eg, to \cite{gadat2018stochastic,gitman2019understanding,gess2026exponential,yuan2016influence,wang2026generalized}.
In particular,
\cite{gadat2018stochastic} analyzes a stochastic heavy-ball scheme for smooth,
coercive, and globally strongly convex objectives and derives mean square bounds
of order $\mathcal{O}(\gamma_n)$ for polynomially decaying learning rates. Moreover, in the quadratic case with a time-varying momentum coefficient, it derives a bound that separates an exponentially decaying initial-condition term from a stochastic term of order $\gamma_n$.
\cite{gitman2019understanding} uses the
quasi-hyperbolic momentum framework to study asymptotic convergence
conditions, stability regions, and stationary distributions for several
stochastic momentum methods. %
In \cite{gess2026exponential}, the authors use a related Lyapunov argument to prove exponential convergence for momentum \SGD\ under a local Polyak--Lojasiewicz condition and a noise assumption  under which the variance vanishes with the objective value.
Recently, \cite{wang2026generalized} introduced a generalized \SGD\ framework with momentum and derived
objective-gap and stationarity guarantees for smooth convex and nonconvex
objectives.

For a broader overview, we refer to \cite{garrigos2023handbook}. In particular, for smooth convex finite-sum objectives, \cite[Theorem~7.4]{garrigos2023handbook} establishes a bound for stochastic momentum with prescribed time-varying momentum and step size parameters.
In the small constant-step regime,
\cite{yuan2016influence} shows that stochastic heavy-ball and Nesterov methods are approximately equivalent to standard \SGD\ with
a rescaled step size. 
Similarly, for fixed-parameter momentum \SGD, \cite{liu2020improved} establishes convergence bounds of the same order as those for \SGD\ for smooth strongly convex and nonconvex objectives and also analyzes multistage parameter schedules.

For a unified analysis of \SGD, stochastic heavy-ball and Nesterov momentum in convex and nonconvex settings, we refer to \cite{liu2024almost,yang2016unified}.

Mini-batch stochastic momentum methods are studied, \eg, in
\cite{tang2023acceleration,bollapragada2025fast}.
\cite{tang2023acceleration} derives finite-sample and averaging results with
explicit batch-size dependence, whereas
\cite{bollapragada2025fast} proves accelerated linear convergence for
sufficiently large mini-batches on quadratic least-squares problems.

\subsection{Structure of this article}

The remainder of this work is organized as follows.
In \cref{section:2} we first present the framework that we employ to study the momentum \SGD\ method (\cref{section:SGD}), we derive upper bounds for the second moment of the auxiliary momentum process, for the correlation between the momentum optimization process and the auxiliary momentum process, and for the mean square error of the momentum \SGD\ method (\cref{section:estimates}), and we employ these upper bounds to derive a Lyapunov recursion for the momentum \SGD\ method (\cref{section:Lyap}); see \cref{prop:Vn:recursionN:eq1} in \cref{prop:Vn:recursionN} for details.

In \cref{section:3} we apply this Lyapunov recursion to different classes of step size sequences. We prove rates of convergence for a class of sublinearly decaying step size sequences and show three regimes for inverse-linear step sizes. Finally,
we prove mean square convergence
for general nonincreasing learning rate sequences which converge to zero and whose sum
diverges.

\section{Mini-batch SGD with momentum: framework and Lyapunov analysis}\label{section:2}
In \cref{section:SGD}, we introduce the key mathematical objects used throughout the paper. 
We specifically introduce the \SGD\ with momentum algorithm along with additional objects required in later sections.
In \cref{section:estimates} we derive separate estimates for
the momentum second moment, the optimization process--momentum correlation, and the mean square error. 
In \cref{section:Lyap} we combine these estimates in an augmented Lyapunov functional
and derive the recursive inequality which will be used to obtain the strong
error bounds.

\subsection{SGD with momentum framework}\label{section:SGD}
In this subsection, we introduce the key mathematical concepts used throughout the paper in \cref{setting2} and discuss their interpretations in \cref{setting2:remark}.

\begin{samepage}
\begin{tcolorbox}[colback=white!95!gray,
                  colframe=black,
                  boxrule=0.5pt,
                  sharp corners,
                  enhanced,
                  breakable,
                 ]
\begin{setting}\label{setting2}
    Let $ ( \Omega, \mathcal{F}, \Pp)$ be a probability space, 
    let $\fd, \dimm, J \in \N$,
    let $\Inn_{n,j} \colon \Omega \to \R^\dimm$, $n \in \N_0$, $j \in \{1,\ldots,J\}$, be \iid\ random variables, 
    let $\Gbatch \colon \R^\fd \times \R^\dimm \to \R^\fd $ be measurable,
    assume for all  $\theta \in \R^\fd$ that $\E\bigl[\|\Gbatch(\theta,\Inn_{0,1})\|\bigr]< \infty$, 
    let $\f \colon \R^\fd \to \R^\fd$ satisfy for all $\theta \in \R^\fd$ that $\f(\theta) = \E[\Gbatch(\theta,\Inn_{0,1})]$,
    let $\alpha \in [0,1)$, $ \nu \in \R^\fd$,
    let $(\gamma_n)_{n \in \N_0}\subseteq (0,\infty)$,
    let 
    $ \Theta \colon \N_0 \times \Omega \to \R^\fd $ and $ \m \colon \N_0 \times \Omega \to \R^\fd $
    be stochastic processes which
    satisfy for all $ n \in \N $ that 
    \begin{equation}
    \m_0=\nu, \qquad \Theta_n 
      = \Theta_{ n - 1 } + \gamma_n \m_n,
      \end{equation}
    \begin{equation}
      \andq \m_n=\alpha \m_{n-1}+(1-\alpha)  \left[\frac{1}{J} \sum_{j=1}^{J} \Gbatch(\Theta_{n-1}, \Inn_{n,j})\right],
    \end{equation}
    assume that $\Theta_0$ and $( \Inn_{ n,j } )_{ n \in \N_0, j \in\{1,\ldots,J\} }$ are independent, assume that $\E[\|\Theta_0\|^2]<\infty$,
    for every $n \in \N_0$ let $\cF_n = \sigma\big(\{\Theta_k \colon k \in \{0,1,\dots,n\}\}\big)$,
    assume for all $n \in \N$, $j \in\{1,\ldots,J\} $ that $\Inn_{n,j} $ is independent of $\cF_{n-1}$
    and let $\const \in [0,\infty]$ satisfy $\const \geq \sup_{n \in \N_0}(\E[ \|\Gbatch(\Theta_n, \Inn_{0,1})\|^2])^{\nicefrac12}$.
\end{setting}
\end{tcolorbox}
\end{samepage}

\begin{remark}[Explanation for \cref{setting2}]\label{setting2:remark}
    In this remark we provide some intuitive interpretations for the mathematical objects appearing in \cref{setting2}. 
    Roughly speaking, we note that
    \begin{enumerate}[label=(\roman*)]
    \item we think of $\fd$ as the dimension of the underlying \SAP,
    \item we think of $\dimm$ as the dimension of the space in which the data (random variables) of the underlying \SAP\ take values,
    \item we think of $\Gbatch$ as the negative gradient of loss function that we aim to minimize,
    \item we think of $( \Inn_{ n,j } )_{ (n,j) \in \N_0 \times \{1,\dots,J\} }$ as random samples used in the algorithm for training,
    \item we think of $\alpha$ as the damping factor of the \SGD\ with momentum algorithm,
    \item we think of $ \nu$ as the initialization of the momentum,
    \item we think of $(\gamma_n)_{n\in \N_0}$ as the sequence of step sizes of the \SGD\ with momentum algorithm,
    \item we think of $(\Theta_n)_{n\in \N_0}$ as the \SGD\ with momentum training process, and
    \item we think of $\f$ as the \SGD\ with momentum vector field.
    \end{enumerate}
\end{remark}

\subsection{Moment, correlation, and error estimates}\label{section:estimates}
    In this subsection in \cref{lem:momentum-second-moment}, \cref{lem:r-recursion2}, and \cref{prop:err2} we derive bounds 
    for the second moment of the auxiliary momentum process, for the correlation between the momentum optimization process and the auxiliary momentum process, and for the mean square error.
    
\begin{athm}{lemma}{lem:momentum-second-moment}
    Assume \cref{setting2}, assume $\const<\infty$, and let $c \in(0, \infty)$, $\vartheta \in \R^\fd$ satisfy for all $\theta \in \R^\fd$ that
    \begin{equation}
        \langle\theta-\vartheta, \E[\Gbatch(\theta, \Inn_{1,1})]\rangle \leq -c \max\{\|\theta-\vartheta\|^2,\|\E[\Gbatch(\theta,\Inn_{1,1})]\|^2\}.
    \end{equation}
    Then it holds for all $n\in\N$ that
    \[
        \E[\|\m_n\|^2]
        \leq
        \alpha\E[\|\m_{n-1}\|^2]
        +
        \frac{(1-\alpha)(J-1+\alpha)}{Jc^2}
        \E[\|\Theta_{n-1}-\vartheta\|^2]
        +
        \frac{(1-\alpha)^2\const^2}{J}.
    \]
\end{athm}
\begin{aproof}
Throughout this proof for every $n\in\N$ let $\varepsilon_n \colon \Omega \to \R^{\fd}$ satisfy
\begin{equation}\llabel{lem:momentum-second-moment:epsilon-def}
\varepsilon_n
=
\frac{1}{J}\smallsum_{j=1}^J
\Gbatch(\Theta_{n-1},\Inn_{n,j})
-
\f(\Theta_{n-1}).
\end{equation}
\argument{\lref{lem:momentum-second-moment:epsilon-def}}{%
that for all $n\in\N$ it holds that
\begin{equation}\llabel{lem:momentum-second-moment:m-dec}
\begin{split}
    \m_n
    &=
    \alpha\m_{n-1}
    +
    \frac{1-\alpha}{J}
    \smallsum_{j=1}^J
    \Gbatch(\Theta_{n-1},\Inn_{n,j})
    \\
    &=
    \alpha\m_{n-1}
    +(1-\alpha)\f(\Theta_{n-1})
    +(1-\alpha)\varepsilon_n.
\end{split}
\end{equation}}
\argument{the fact that $\Theta_{n-1}$ is $\cF_{n-1}$-measurable;
the independence of $\Inn_{n,1},\ldots,\Inn_{n,J}$ from $\cF_{n-1}$}{%
that for all $n\in\N$ it holds that
\begin{equation}\llabel{lem:momentum-second-moment:epsilon-centered}
\begin{split}
    \E[\varepsilon_n\mid\cF_{n-1}]
    &=
    \frac{1}{J}\smallsum_{j=1}^J
    \E[
    \Gbatch(\Theta_{n-1},\Inn_{n,j})
    \mid\cF_{n-1}]
    -
    \f(\Theta_{n-1})
    \\
    &=
    \frac{1}{J}\smallsum_{j=1}^J
    \f(\Theta_{n-1})
    -
    \f(\Theta_{n-1})
    \\
    &=0.
\end{split}
\end{equation}}
\argument{the fact that $\m_{n-1}$ and $\f(\Theta_{n-1})$ are
$\cF_{n-1}$-measurable;
\lref{lem:momentum-second-moment:epsilon-centered};
the tower property}{%
that for all $n\in\N$ it holds that
\begin{equation}\llabel{lem:momentum-second-moment:cross}
\begin{split}
    &\E\!\left[
    \left\langle
    \alpha\m_{n-1}+(1-\alpha)\f(\Theta_{n-1}),
    \varepsilon_n
    \right\rangle
    \right]
    \\
    &=
    \E\!\left[
    \E\!\left[
    \left\langle
    \alpha\m_{n-1}+(1-\alpha)\f(\Theta_{n-1}),
    \varepsilon_n
    \right\rangle
    \,\middle|\,
    \cF_{n-1}
    \right]
    \right]
    \\
    &=
    \E\!\left[
    \left\langle
    \alpha\m_{n-1}+(1-\alpha)\f(\Theta_{n-1}),
    \E[\varepsilon_n\mid\cF_{n-1}]
    \right\rangle
    \right]
    \\
    &=0.
\end{split}
\end{equation}}
\argument{\lref{lem:momentum-second-moment:m-dec};
\lref{lem:momentum-second-moment:cross}}{%
that for all $n\in\N$ it holds that
\begin{equation}\llabel{lem:momentum-second-moment:dec}
\begin{split}
    \E[\|\m_n\|^2]
    &=
    \E\!\left[
    \|
    \alpha\m_{n-1}
    +(1-\alpha)\f(\Theta_{n-1})
    +(1-\alpha)\varepsilon_n
    \|^2
    \right]
    \\
    &=
    \E\!\left[
    \|\alpha\m_{n-1}
    +(1-\alpha)\f(\Theta_{n-1})\|^2
    \right]
    \\
    &\quad
    +2(1-\alpha)
    \E\!\left[
    \left\langle
    \alpha\m_{n-1}
    +(1-\alpha)\f(\Theta_{n-1}),
    \varepsilon_n
    \right\rangle
    \right]
    \\
    &\quad
    +(1-\alpha)^2\E[\|\varepsilon_n\|^2]
    \\
    &=
    \E\!\left[
    \|\alpha\m_{n-1}
    +(1-\alpha)\f(\Theta_{n-1})\|^2
    \right]
    +(1-\alpha)^2\E[\|\varepsilon_n\|^2].
\end{split}
\end{equation}}
\argument{the convexity of the squared norm;
\lref{lem:momentum-second-moment:dec}}{%
that for all $n\in\N$ it holds that
\begin{equation}\llabel{lem:momentum-second-moment:conv}
\begin{split}
    \E[\|\m_n\|^2]
    &\leq
    \alpha\E[\|\m_{n-1}\|^2]
    +(1-\alpha)\E[\|\f(\Theta_{n-1})\|^2]
    \\
    &\quad
    +(1-\alpha)^2\E[\|\varepsilon_n\|^2].
\end{split}
\end{equation}}
\argument{\lref{lem:momentum-second-moment:epsilon-def};
the fact that
$\E[\Gbatch(\Theta_{n-1},\Inn_{n,1})\mid\cF_{n-1}]
=\f(\Theta_{n-1})$;
the tower property}{%
that for all $n\in\N$ it holds that
\begin{equation}\llabel{lem:momentum-second-moment:epsilon-var}
\begin{split}
    \E[\|\varepsilon_n\|^2]
    &=
    \frac{1}{J^2}
    \E\!\left[
    \left\|
    \smallsum_{j=1}^J
    \left(
    \Gbatch(\Theta_{n-1},\Inn_{n,j})
    -
    \f(\Theta_{n-1})
    \right)
    \right\|^2
    \right]
    \\
    &=
    \frac{1}{J^2}
    \smallsum_{j=1}^J
    \E\!\left[
    \|
    \Gbatch(\Theta_{n-1},\Inn_{n,j})
    -
    \f(\Theta_{n-1})
    \|^2
    \right]
    \\
    &=
    \frac{1}{J}
    \E\!\left[
    \|
    \Gbatch(\Theta_{n-1},\Inn_{n,1})
    -
    \f(\Theta_{n-1})
    \|^2
    \right]
    \\
    &=
    \frac{1}{J}
    \bigg(
    \E[\|\Gbatch(\Theta_{n-1},\Inn_{n,1})\|^2]
    \\
    &\qquad
    -2\E\!\left[
    \left\langle
    \Gbatch(\Theta_{n-1},\Inn_{n,1}),
    \f(\Theta_{n-1})
    \right\rangle
    \right]
    +
    \E[\|\f(\Theta_{n-1})\|^2]
    \bigg)
    \\
    &=
    \frac{1}{J}
    \bigg(
    \E[\|\Gbatch(\Theta_{n-1},\Inn_{n,1})\|^2]
    \\
    &\qquad
    -2\E\!\left[
    \left\langle
    \E[
    \Gbatch(\Theta_{n-1},\Inn_{n,1})
    \mid\cF_{n-1}],
    \f(\Theta_{n-1})
    \right\rangle
    \right]
    +
    \E[\|\f(\Theta_{n-1})\|^2]
    \bigg)
    \\
    &=
    \frac{1}{J}
    \left(
    \E[\|\Gbatch(\Theta_{n-1},\Inn_{n,1})\|^2]
    -
    \E[\|\f(\Theta_{n-1})\|^2]
    \right)
    \\
    &\leq
    \frac{1}{J}
    \left(
    \const^2
    -
    \E[\|\f(\Theta_{n-1})\|^2]
    \right).
\end{split}
\end{equation}}

\argument{\lref{lem:momentum-second-moment:conv};
\lref{lem:momentum-second-moment:epsilon-var}}{%
that for all $n\in\N$ it holds that
\begin{equation}\llabel{lem:momentum-second-moment:f-bound}
\begin{split}
    \E[\|\m_n\|^2]
    &\leq
    \alpha\E[\|\m_{n-1}\|^2]
    +(1-\alpha)\E[\|\f(\Theta_{n-1})\|^2]
    \\
    &\quad
    +
    \frac{(1-\alpha)^2}{J}
    \left(
    \const^2
    -
    \E[\|\f(\Theta_{n-1})\|^2]
    \right)
    \\
    &=
    \alpha\E[\|\m_{n-1}\|^2]
    +
    \left(
    1-\alpha-\frac{(1-\alpha)^2}{J}
    \right)
    \E[\|\f(\Theta_{n-1})\|^2]
    \\
    &\quad
    +
    \frac{(1-\alpha)^2\const^2}{J}
    \\
    &=
    \alpha\E[\|\m_{n-1}\|^2]
    +
    \frac{J(1-\alpha)-(1-\alpha)^2}{J}
    \E[\|\f(\Theta_{n-1})\|^2]
    \\
    &\quad
    +
    \frac{(1-\alpha)^2\const^2}{J}
    \\
    &=
    \alpha\E[\|\m_{n-1}\|^2]
    +
    \frac{(1-\alpha)(J-1+\alpha)}{J}
    \E[\|\f(\Theta_{n-1})\|^2]
    \\
    &\quad
    +
    \frac{(1-\alpha)^2\const^2}{J}.
\end{split}
\end{equation}}
\argument{the assumption that for all $\theta\in\R^\fd$ it holds that
\begin{equation}
    \langle\theta-\vartheta,\f(\theta)\rangle
    \leq
    -c\max\{
    \|\theta-\vartheta\|^2,
    \|\f(\theta)\|^2
    \}
\end{equation};
the Cauchy--Schwarz inequality}{%
that for all $\theta\in\R^\fd$ it holds that
\begin{equation}\llabel{lem:momentum-second-moment:coer}
\begin{split}
    c\|\f(\theta)\|^2
    &\leq
    -\langle\theta-\vartheta,\f(\theta)\rangle
    \\
    &\leq
    \|\theta-\vartheta\|\|\f(\theta)\|.
\end{split}
\end{equation}}
\argument{\lref{lem:momentum-second-moment:f-bound};
\lref{lem:momentum-second-moment:coer}}{%
that for all $n\in\N$ it holds that
\begin{equation}
\begin{split}
    \E[\|\m_n\|^2]
    &\leq
    \alpha\E[\|\m_{n-1}\|^2]
    +
    \frac{(1-\alpha)(J-1+\alpha)}{Jc^2}
    \E[\|\Theta_{n-1}-\vartheta\|^2]
    +
    \frac{(1-\alpha)^2\const^2}{J}.
\end{split}
\end{equation}}
\end{aproof}

\begin{athm}{lemma}{lem:r-recursion2}
    Assume \cref{setting2}, assume $\const < \infty$, and let $c \in(0, \infty)$, $\vartheta \in \R^\fd$ satisfy for all $\theta \in \R^\fd$ that
    \begin{equation}
        \langle\theta-\vartheta, \E[\Gbatch(\theta, \Inn_{1,1})]\rangle \leq -c \max\{\|\theta-\vartheta\|^2,\|\E[\Gbatch(\theta,\Inn_{1,1})]\|^2\}.
    \end{equation}
    Then it holds for all $n \in \N$ that
    \begin{equation}
         \E[\langle \Theta_{n-1} - \vartheta, \m_n\rangle] \leq \alpha \E[\langle \Theta_{n-1} - \vartheta, \m_{n-1}\rangle] - (1-\alpha)c  \E[\|\Theta_{n-1} - \vartheta\|^2]
    \end{equation}
    and
    \begin{equation}
        \max\{\E[\langle \Theta_n - \vartheta, \m_n\rangle],0\} \leq \alpha \max\{\E[\langle \Theta_{n-1} - \vartheta, \m_{n-1}\rangle],0\}   + \gamma_n \E[\|\m_n\|^2].
    \end{equation}
\end{athm}
\begin{aproof}
    \argument{ $\Theta_n - \vartheta = (\Theta_{n-1}-\vartheta) + \gamma_n \m_n$}{%
    that for all $n \in \N$ it holds that
    \begin{equation}\llabel{lem:r-recursion:dec}
         \E[\langle \Theta_n - \vartheta, \m_n\rangle]  = \E[\langle \Theta_{n-1}-\vartheta, \m_n\rangle] + \gamma_n \E[\|\m_n\|^2].
    \end{equation}}
    \argument{the fact that $\Theta_{n-1}-\vartheta$ is $\cF_{n-1}$-measurable;
    the fact that $\E[\m_n\mid\cF_{n-1}] = \alpha \m_{n-1} + (1-\alpha)\E[(\nicefrac1J ) \sum_{j=1}^J\Gbatch(\Theta_{n-1},\Inn_{n,j})\mid \mathcal F_{n-1}] = \alpha\m_{n-1} + (1-\alpha)\f(\Theta_{n-1})$;
    the tower property}{that for all $n \in \N$ it holds that
    \begin{equation}\llabel{lem:r-recursion:tow}
    \begin{aligned}
    \E[\langle \Theta_{n-1}-\vartheta, \m_n\rangle]
    &= \E\Big[\E\big[\langle \Theta_{n-1}-\vartheta, \m_n\rangle \,\big|\, \mathcal F_{n-1}\big]\Big] \\
    &= \E\Big[\big\langle \Theta_{n-1}-\vartheta,\E[\m_n \mid \mathcal F_{n-1}]\big\rangle\Big] \\
    &= \E\Big[\big\langle \Theta_{n-1}-\vartheta,\alpha \m_{n-1} + (1-\alpha)\f(\Theta_{n-1})\big\rangle\Big] \\
    &= \alpha \E[\langle \Theta_{n-1}-\vartheta, \m_{n-1} \rangle]
    + (1-\alpha)\E[\langle \Theta_{n-1}-\vartheta, \f(\Theta_{n-1})\rangle].
    \end{aligned}
    \end{equation}}
    \argument{the assumption that for all $\theta \in \R^\fd$ it holds
    that $\langle\theta-\vartheta, \f(\theta)\rangle \leq -c\|\theta-\vartheta\|^2$}{that for all $n \in \N$ it holds that
    \begin{equation}\llabel{lem:r-recursion:coer}
        \E[\langle \Theta_{n-1}-\vartheta, \f(\Theta_{n-1})\rangle] \leq -c \E[\|\Theta_{n-1} - \vartheta\|^2].
    \end{equation}}
    \argument{\lref{lem:r-recursion:tow}; \lref{lem:r-recursion:coer}}{%
    that for all $n \in \N$ it holds that
    \begin{equation}\llabel{lem:r-recursion:1}
         \E[\langle \Theta_{n-1} - \vartheta, \m_n\rangle] \leq \alpha \E[\langle \Theta_{n-1} - \vartheta, \m_{n-1}\rangle] - (1-\alpha)c  \E[\|\Theta_{n-1} - \vartheta\|^2].
    \end{equation}
    }
    \argument{\lref{lem:r-recursion:dec}; \lref{lem:r-recursion:1}}{
    \begin{equation}\llabel{lem:r-recursion:2noplus}
        \E[\langle \Theta_n - \vartheta, \m_n\rangle] \leq \alpha\E[\langle \Theta_{n-1} - \vartheta, \m_{n-1}\rangle] - (1-\alpha) c  \E[\|\Theta_{n-1} - \vartheta\|^2] + \gamma_n \E[\|\m_n\|^2].
    \end{equation}
    }
    \argument{\lref{lem:r-recursion:2noplus}}{
    that for all $n \in \N$ it holds that
    \begin{equation}
    \begin{split}
        \max\{\E[\langle \Theta_n - \vartheta, \m_n\rangle],0\} & \leq \max\{\alpha\E[\langle \Theta_{n-1} - \vartheta, \m_{n-1}\rangle] - (1-\alpha) c  \E[\|\Theta_{n-1} - \vartheta\|^2] \\
        &\quad + \gamma_n \E[\|\m_n\|^2],0\}
        \\
        &\leq \alpha \max\{\E[\langle \Theta_{n-1} - \vartheta, \m_{n-1}\rangle],0\}   + \gamma_n \E[\|\m_n\|^2].   
    \end{split}
    \end{equation}
    }
\end{aproof}

\begin{athm}{lemma}{prop:err2}
    Assume \cref{setting2}, assume $\const < \infty$,  and let $c \in(0, \infty)$, $\vartheta \in \R^\fd$ satisfy for all $\theta \in \R^\fd$ that
    \begin{equation}
        \langle\theta-\vartheta, \E[\Gbatch(\theta, \Inn_{1,1})]\rangle \leq -c \max\{\|\theta-\vartheta\|^2,\|\E[\Gbatch(\theta,\Inn_{1,1})]\|^2\}.
    \end{equation}
    Then it holds for all $n \in \N$ that
        \begin{equation}
        \begin{split}
             \E[\|\Theta_n - \vartheta\|^2] &\leq (1 - 2(1-\alpha)c\gamma_n)\E[\|\Theta_{n-1}- \vartheta\|^2] + 2\alpha\gamma_n \E[\langle \Theta_{n-1} - \vartheta, \m_{n-1}\rangle]\\
             & \quad + \gamma_n^2 \E[\|\m_n\|^2].
        \end{split}
        \end{equation} 
\end{athm}
\begin{aproof}
    \argument{\cref{lem:r-recursion2}}{that for all $n\in\N$ it holds that
    \begin{equation}\llabel{prop2:err:1}
    \begin{split}
         \E[\|\Theta_n - \vartheta\|^2] &=  \E[\|\Theta_{n-1} - \vartheta\|^2] + 2\gamma_n\E[\langle\Theta_{n-1}-\vartheta,\m_n\rangle] + \gamma_n^2\E[\|\m_n\|^2]\\
         & \leq \E[\|\Theta_{n-1} - \vartheta\|^2] + 2\gamma_n \big( \alpha \E[\langle \Theta_{n-1} - \vartheta, \m_{n-1}\rangle] 
         \\
         & \quad - (1-\alpha)c \, \E[\|\Theta_{n-1} - \vartheta\|^2] \big)  + \gamma_n^2\E[\|\m_n\|^2]
         \\
        &= (1 - 2(1-\alpha)c\,\gamma_n) \E[\|\Theta_{n-1} - \vartheta\|^2] + 2\alpha\gamma_n \E[\langle \Theta_{n-1} - \vartheta, \m_{n-1}\rangle] \\
        & \quad +\gamma_n^2\E[\|\m_n\|^2].
    \end{split}
    \end{equation}}
\end{aproof}

\subsection{Augmented Lyapunov recursion}\label{section:Lyap}
We combine the three bounds obtained in \cref{section:estimates} to derive an augmented Lyapunov functional. \cref{prop:Vn} derives a
one-step estimate for this functional, \cref{prop:Vn:recursion} rewrites the estimate as a scalar recursive inequality, and \cref{prop:Vn:recursionN} simplifies the recursion for nonincreasing and sufficiently small learning rates.

\begin{athm}{lemma}{prop:Vn}
    Assume \cref{setting2}, assume $\const < \infty$,  and let $c \in(0, \infty)$, $\vartheta \in \R^\fd$ satisfy for all $\theta \in \R^\fd$ that
    \begin{equation}
        \langle\theta-\vartheta, \E[\Gbatch(\theta, \Inn_{1,1})]\rangle \leq -c \max\{\|\theta-\vartheta\|^2,\|\E[\Gbatch(\theta,\Inn_{1,1})]\|^2\}.
    \end{equation}
    Then it holds for all $n \in \N$ that
        \begin{equation}\label{prop:Vn:rec}
        \begin{split}
             &\E[\|\Theta_n - \vartheta\|^2] + \frac{1+\alpha}{1-\alpha} \gamma_n \max\{\E[\langle \Theta_{n} - \vartheta, \m_{n}\rangle],0\} + \frac{2}{(1-\alpha)^2} \gamma_n^2 \E[\|\m_n\|^2]  \\
             & \leq \left(1-2(1-\alpha)c\gamma_n + \frac{2(2-\alpha)(J-1+\alpha)}{(1-\alpha)Jc^2}\gamma_n^2\right)\E[\|\Theta_{n-1} - \vartheta\|^2]\\
             & \quad +\frac{\alpha(3-\alpha)}{1-\alpha}\gamma_n  \max\{\E[\langle \Theta_{n-1} - \vartheta, \m_{n-1}\rangle],0\} + \frac{2\alpha(2-\alpha)}{(1-\alpha)^2}\gamma_n^2 \E[\|\m_{n-1}\|^2] \\
             & \quad + \frac{2(2-\alpha)\const^2}{J} \gamma_n^2.
        \end{split}
        \end{equation} 
\end{athm}
\begin{aproof}
Throughout this proof let $p,q,H \in(0,\infty)$ satisfy
\begin{equation}
    p=\frac{1+\alpha}{1-\alpha},
    \qquad
    q=\frac{2}{(1-\alpha)^2},
    \qquad
    H=1+p+q.
\end{equation}
\argument{\cref{prop:err2}; \cref{lem:r-recursion2}}{%
that for all $n\in\N$ it holds that
\begin{equation}\llabel{prop:Vn:first}
\begin{split}
    &\E[\|\Theta_n-\vartheta\|^2]
    +p\gamma_n
    \max\{\E[\langle\Theta_n-\vartheta,\m_n\rangle],0\}
    +q\gamma_n^2\E[\|\m_n\|^2]
    \\
    &\leq
    \bigl(1-2(1-\alpha)c\gamma_n\bigr)
    \E[\|\Theta_{n-1}-\vartheta\|^2]
    +2\alpha\gamma_n
    \E[
      \langle\Theta_{n-1}-\vartheta,\m_{n-1}\rangle
    ]
    +\gamma_n^2\E[\|\m_n\|^2]
    \\
    &\quad
    +p\gamma_n
    \bigl[
      \alpha
      \max\{\E[
        \langle\Theta_{n-1}-\vartheta,\m_{n-1}\rangle
      ],0\}
      +\gamma_n\E[\|\m_n\|^2]
    \bigr]
    +q\gamma_n^2\E[\|\m_n\|^2]
    \\
    & \leq
    \bigl(
      1-2(1-\alpha)c\gamma_n
    \bigr)
    \E[\|\Theta_{n-1}-\vartheta\|^2]
    \\
    &\quad
    +\alpha(2+p)\gamma_n
    \max\{\E[
      \langle\Theta_{n-1}-\vartheta,\m_{n-1}\rangle
    ],0\}
    +H\gamma_n^2\E[\|\m_n\|^2].
\end{split}
\end{equation}}
\argument{\lref{prop:Vn:first};
\cref{lem:momentum-second-moment}}{%
that for all $n\in\N$ it holds that
\begin{equation}\llabel{prop:Vn:second}
\begin{split}
    &\E[\|\Theta_n-\vartheta\|^2]
    +p\gamma_n
    \max\{\E[\langle\Theta_n-\vartheta,\m_n\rangle],0\}
    +q\gamma_n^2\E[\|\m_n\|^2]
    \\
    &\leq
    \left[
      1-2(1-\alpha)c\gamma_n
      +
      H
      \frac{(1-\alpha)(J-1+\alpha)}{Jc^2}
      \gamma_n^2
    \right]
    \E[\|\Theta_{n-1}-\vartheta\|^2]
    \\
    &\quad
    +\alpha(2+p)\gamma_n
    \max\{\E[
      \langle\Theta_{n-1}-\vartheta,\m_{n-1}\rangle
    ],0\}
    \\
    &\quad
    +\alpha H\gamma_n^2
    \E[\|\m_{n-1}\|^2]
    +
    H\frac{(1-\alpha)^2\const^2}{J}\gamma_n^2.
\end{split}
\end{equation}}
\argument{the fact that
\begin{equation}
    2+p=\frac{3-\alpha}{1-\alpha},
    \qquad
    H=1+p+q
    =
    \frac{2(2-\alpha)}{(1-\alpha)^2}
\end{equation};
\lref{prop:Vn:second}}{%
that for all $n\in\N$ it holds that
\begin{equation}
\begin{split}
    &\E[\|\Theta_n-\vartheta\|^2]
    +\frac{1+\alpha}{1-\alpha}\gamma_n
    \max\{\E[\langle\Theta_n-\vartheta,\m_n\rangle],0\}
    +\frac{2}{(1-\alpha)^2}\gamma_n^2
    \E[\|\m_n\|^2]
    \\
    &\leq
    \left[
      1-2(1-\alpha)c\gamma_n
      +
      \frac{
        2(2-\alpha)(J-1+\alpha)
      }{
        (1-\alpha)Jc^2
      }
      \gamma_n^2
    \right]
    \E[\|\Theta_{n-1}-\vartheta\|^2]
    \\
    &\quad
    +
    \frac{\alpha(3-\alpha)}{1-\alpha}
    \gamma_n
    \max\{\E[
      \langle\Theta_{n-1}-\vartheta,\m_{n-1}\rangle
    ],0\}
    \\
    &\quad
    +
    \frac{2\alpha(2-\alpha)}{(1-\alpha)^2}
    \gamma_n^2
    \E[\|\m_{n-1}\|^2]
    +
    \frac{2(2-\alpha)\const^2}{J}\gamma_n^2.
\end{split}
\end{equation}}
\end{aproof}
\begin{remark}[Choice of the Lyapunov coefficients]
In this remark we provide some intuition for the Lyapunov coefficients appearing in \cref{prop:Vn}.
The coefficients
\begin{equation}
    \frac{1+\alpha}{1-\alpha}
    \qquad\text{and}\qquad
    \frac{2}{(1-\alpha)^2}
\end{equation}
of the second and third terms in \cref{prop:Vn:rec}
are chosen so that the correlation and momentum terms arising can be absorbed into the corresponding terms of
$V_{n-1}$. Indeed, it holds that
\begin{equation}
    \frac{\alpha(3-\alpha)}{1-\alpha}= \frac{\alpha(3-\alpha)}{1+\alpha}\frac{1+\alpha}{1-\alpha} < \frac{1+\alpha}{1-\alpha}
\end{equation}
\begin{equation}
    \andq \frac{2\alpha(2-\alpha)}{(1-\alpha)^2} = \alpha(2-\alpha) \frac{2}{(1-\alpha)^2}
    <\frac{2}{(1-\alpha)^2}.
\end{equation}
Consequently, for nonincreasing step sizes and sufficiently small
$\gamma_n$, both terms can be bounded using a common coefficient, which yields a scalar recursion. Note that these choices are not unique, but lead to
particularly simple explicit constants.
\end{remark}

\begin{athm}{lemma}{prop:Vn:recursion}
    Assume \cref{setting2}, assume $\const < \infty$, let $c \in(0, \infty)$, $\vartheta \in \R^\fd$ satisfy for all $\theta \in \R^\fd$ that
    \begin{equation}
        \langle\theta-\vartheta, \E[\Gbatch(\theta, \Inn_{1,1})]\rangle \leq -c \max\{\|\theta-\vartheta\|^2,\|\E[\Gbatch(\theta,\Inn_{1,1})]\|^2\},
    \end{equation}
    for every $n \in \N_0$ let $V_n \in [0,\infty)$ satisfy \begin{equation}
        V_n = \E[\|\Theta_n - \vartheta\|^2] + \frac{1+\alpha}{1-\alpha} \gamma_n \max\{\E[\langle \Theta_{n} - \vartheta, \m_{n}\rangle],0\} + \frac{2}{(1-\alpha)^2} \gamma_n^2 \E[\|\m_n\|^2],
    \end{equation}
    and for every $n \in \N$ let $\beta_n \in \R$ satisfy $\beta_n = 1-2(1-\alpha)c\gamma_n + \frac{2(2-\alpha)(J-1+\alpha)}{(1-\alpha)Jc^2}\gamma_n^2$.
    Then it holds for all $n \in \N$ that
        \begin{equation}
        \begin{split}
             & V_n \leq \max\left\{ \beta_n,  \frac{\alpha(3-\alpha)\gamma_n}{(1+\alpha)\gamma_{n-1}}, \frac{\alpha(2-\alpha)\gamma_n^2}{\gamma_{n-1}^2}  \right\} V_{n-1}+ \frac{2(2-\alpha)\const^2}{J} \gamma_n^2
        \end{split}
        \end{equation} 
\end{athm}
\begin{aproof}
Throughout this proof for every $n\in\N_0$ let $e_n,r_n^+,s_n \in [0,\infty)$ satisfy
\begin{equation}
    e_n=\E[\|\Theta_n-\vartheta\|^2],
\end{equation}
\begin{equation}
    r_n^+
    =
    \max\{
      \E[\langle\Theta_n-\vartheta,\m_n\rangle],
      0
    \},
\end{equation}
and
\begin{equation}
    s_n=\E[\|\m_n\|^2]
\end{equation}
and let $p,q,H \in(0,\infty)$ satisfy
\begin{equation}
    p=\frac{1+\alpha}{1-\alpha},
    \qquad
    q=\frac{2}{(1-\alpha)^2},
    \qquad
    H=1+p+q.
\end{equation}
\argument{\cref{prop:Vn}}{%
that for all $n\in\N$ it holds that
\begin{equation}
\begin{split}\llabel{prop:Vn:recursion:first}
    V_n
    &=
    e_n+p\gamma_nr_n^++q\gamma_n^2s_n
    \\
    &\leq
    \left[
      1-2(1-\alpha)c\gamma_n
      +
      H
      \frac{(1-\alpha)(J-1+\alpha)}{Jc^2}
      \gamma_n^2
    \right]e_{n-1}
    \\
    &\quad
    +\alpha(2+p)\gamma_nr_{n-1}^+
    +\alpha H\gamma_n^2s_{n-1}
    +
    H\frac{(1-\alpha)^2\const^2}{J}\gamma_n^2.
\end{split}
\end{equation}}
\argument{the definition of $\beta_n$;
the fact that
\begin{equation}
  2+p=\frac{3-\alpha}{1-\alpha},
  \qquad
  H=\frac{2(2-\alpha)}{(1-\alpha)^2}
\end{equation};
\lref{prop:Vn:recursion:first}}{%
that for all $n\in\N$ it holds that
\begin{equation}\llabel{prop:Vn:recursion:third}
\begin{split}
    V_n
    &\leq
    \beta_ne_{n-1}
    +
    \frac{\alpha(3-\alpha)}{1-\alpha}
    \gamma_nr_{n-1}^+
    \\
    &\quad
    +
    \frac{2\alpha(2-\alpha)}{(1-\alpha)^2}
    \gamma_n^2s_{n-1}
    +
    \frac{2(2-\alpha)\const^2}{J}\gamma_n^2.
\end{split}
\end{equation}}
\argument{the fact that for all $n\in\N$ it holds that
\begin{equation}
\begin{aligned}
    \frac{\alpha(3-\alpha)}{1-\alpha}
    \gamma_nr_{n-1}^+
    &=
    \frac{
      \alpha(3-\alpha)\gamma_n
    }{
      (1+\alpha)\gamma_{n-1}
    }
    \left[
      \frac{1+\alpha}{1-\alpha}
      \gamma_{n-1}r_{n-1}^+
    \right],
    \\
    \frac{2\alpha(2-\alpha)}{(1-\alpha)^2}
    \gamma_n^2s_{n-1}
    &=
    \frac{
      \alpha(2-\alpha)\gamma_n^2
    }{
      \gamma_{n-1}^2
    }
    \left[
      \frac{2}{(1-\alpha)^2}
      \gamma_{n-1}^2s_{n-1}
    \right]
\end{aligned}
\end{equation};
the fact that for all $n\in\N$ it holds $e_{n-1}\geq0$, $r_{n-1}^+\geq0$ and $s_{n-1}\geq0$;
\lref{prop:Vn:recursion:third}}{%
that for all $n\in\N$ it holds that
\begin{equation}
\begin{split}
    V_n
    &\leq
    \max\left\{
      \beta_n,
      \frac{
        \alpha(3-\alpha)\gamma_n
      }{
        (1+\alpha)\gamma_{n-1}
      },
      \frac{
        \alpha(2-\alpha)\gamma_n^2
      }{
        \gamma_{n-1}^2
      }
    \right\}
    V_{n-1}
    +
    \frac{2(2-\alpha)\const^2}{J}\gamma_n^2.
\end{split}
\end{equation}}
\end{aproof}

\begin{athm}{cor}{prop:Vn:recursionN}
    Assume \cref{setting2}, assume $\const < \infty$, let $c \in(0, \infty)$, $\vartheta \in \R^\fd$ satisfy for all $\theta \in \R^\fd$ that
    \begin{equation}
        \langle\theta-\vartheta, \E[\Gbatch(\theta, \Inn_{1,1})]\rangle \leq -c \max\{\|\theta-\vartheta\|^2,\|\E[\Gbatch(\theta,\Inn_{1,1})]\|^2\},
    \end{equation}
    assume for all $n\in\N_0$ that $\gamma_{n+1}\leq \gamma_{n}$,
    for every $n \in \N_0$ let $V_n \in [0,\infty)$ satisfy \begin{equation}
        V_n = \E[\|\Theta_n - \vartheta\|^2] + \frac{1+\alpha}{1-\alpha} \gamma_n \max\{\E[\langle \Theta_{n} - \vartheta, \m_{n}\rangle],0\} + \frac{2}{(1-\alpha)^2} \gamma_n^2 \E[\|\m_n\|^2],
    \end{equation}
    and for every $n \in \N$ let $\beta_n \in \R$ satisfy $\beta_n = 1-2(1-\alpha)c\gamma_n + \frac{2(2-\alpha)(J-1+\alpha)}{(1-\alpha)Jc^2}\gamma_n^2$.
    Then 
    for all $n \in \N$ with $\gamma_n\leq\frac{1-\alpha}{2(1+\alpha)c}$ it holds that
        \begin{equation}\label{prop:Vn:recursionN:eq1}
        \begin{split}
             & V_n \leq \beta_n V_{n-1}+ \frac{2(2-\alpha)\const^2}{J} \gamma_n^2.
        \end{split}
        \end{equation} 
\end{athm}
\begin{aproof}
Throughout this proof let $\eta_1,\eta_2 \in[0,\infty)$ satisfy
\begin{equation}
     \eta_1
    =
    \frac{\alpha(3-\alpha)}{1+\alpha}
\end{equation}
and
\begin{equation}
    \eta_2
    =
    \alpha(2-\alpha).
\end{equation}
\argument{the fact that $\alpha\in[0,1)$}{%
that
\begin{equation}\llabel{prop:Vn:recursionN:eta-order}
\begin{split}
    \eta_1-\eta_2
    &=
    \frac{\alpha(1-\alpha)^2}{1+\alpha}
    \geq0
\end{split}
\end{equation}
and
\begin{equation}\llabel{prop:Vn:recursionN:eta-less-one}
\begin{split}
    1-\eta_1
    &=
    \frac{(1-\alpha)^2}{1+\alpha}
    >0.
\end{split}
\end{equation}}
\argument{
\lref{prop:Vn:recursionN:eta-less-one}}{%
that for all $n\in\N$ with $\gamma_n
    \leq
    \frac{1-\alpha}{2(1+\alpha)c}$ it holds that 
\begin{equation}\llabel{prop:Vn:recursionN:beta-lower}
\begin{split}
    \beta_n
    &=
    1-2(1-\alpha)c\gamma_n
    +
    \frac{
      2(2-\alpha)(J-1+\alpha)
    }{
      (1-\alpha)Jc^2
    }
    \gamma_n^2
    \\
    &\geq
    1-2(1-\alpha)c\gamma_n
    \\
    &\geq
    1-\frac{(1-\alpha)^2}{1+\alpha}
    \\
    &=
    \frac{\alpha(3-\alpha)}{1+\alpha}
    =
    \eta_1
    \geq
    \eta_2.
\end{split}
\end{equation}}
\argument{the  assumption that for all $n\in\N$ it holds that $\gamma_{n+1}\leq\gamma_{n}$;
\lref{prop:Vn:recursionN:beta-lower}}{%
that for all $n\in\N$ with $\gamma_n
    \leq
    \frac{1-\alpha}{2(1+\alpha)c}$ it holds that 
\begin{equation}\llabel{prop:Vn:recursionN:max}
\begin{split}
    &\max\left\{
      \beta_n,
      \frac{
        \alpha(3-\alpha)\gamma_n
      }{
        (1+\alpha)\gamma_{n-1}
      },
      \frac{
        \alpha(2-\alpha)\gamma_n^2
      }{
        \gamma_{n-1}^2
      }
    \right\}
    \\
    &\leq
    \max\{\beta_n,\eta_1,\eta_2\}
    =
    \beta_n.
\end{split}
\end{equation}}
\argument{\cref{prop:Vn:recursion};
\lref{prop:Vn:recursionN:max}}{%
that for all $n\in\N$ with $\gamma_n
    \leq
    \frac{1-\alpha}{2(1+\alpha)c}$ it holds that 
\begin{equation}
    V_n
    \leq
    \beta_nV_{n-1}
    +
    \frac{2(2-\alpha)\const^2}{J}\gamma_n^2.
\end{equation}}
\end{aproof}

\section{Strong error bounds and convergence rates} \label{section:3}
In this section 
we derive the main estimate in \cref{theo:Vn}. This estimate is subsequently used in \cref{cor:main-result} to prove \cref{theorem1}.
In \cref{section:rates} we consider polynomially decaying step sizes
while in \cref{section:rates2} we treat the inverse-linear case and obtain
three different regimes. 
In \cref{section:generalgamma} we demonstrate the convergence of the \SGD\ with momentum algorithm to the unique zero of $\f$ with a general sequence of step sizes. The uniqueness of the zero is guaranteed by the coercivity assumption (cf., \eg, \cite[Lemma 2.1]{jentzen2021strong}).

\subsection{Sublinear step sizes}\label{section:rates}

We analyze the rate of convergence of the \SGD\ with momentum algorithm to the zero of the function $\f$ for sublinear step size sequences. Our approach involves expanding the recursive inequality from \cref{prop:Vn:recursionN} by iteratively substituting the recursive terms to derive a closed-form expression. 
This allows us to calculate the convergence rates for each resulting term using the elementary result in \cref{exp}.
In \cref{theo:Vn} we establish the convergence rate.
Only for completeness we include here a proof for \cref{exp}.
\begin{athm}{lemma}{exp}
    Let $\lambda\in (0,1)$, $a,\cgamma\in (0,\infty)$ and let $(\gamma_n)_{n \in \N}\subseteq (0,\infty)$ satisfy for all $n\in\N$ that $\gamma_n=\cgamma n^{-\lambda}$.
    Then 
    \begin{enumerate}[label=(\roman*)]
        \item \label{exp:1} for all $m,n \in \N$ with $(a\cgamma)^{\nicefrac1\lambda}\leq m\leq n$ it holds that
    \begin{equation}
        \prod_{k=m}^n(1-a\gamma_k)\leq \exp\left(-\frac{a\cgamma}{1-\lambda}[(n+1)^{1-\lambda}-m^{1-\lambda}]\right)
    \end{equation}
    and
    \item \label{exp:2} for all $m,n \in \N$ with $\max \left\{(a\cgamma)^{\nicefrac{1}{\lambda}}, \left(\frac{4\lambda}{a\cgamma}\right)^{\nicefrac{1}{1-\lambda}}\right\}\leq m\leq n$ it holds that
    \begin{equation}
        \sum_{k=m}^n\gamma_k^2 \prod_{j=k+1}^n (1-a\gamma_j) \leq \frac{2\cgamma}{a}n^{-\lambda}.
    \end{equation}
    \end{enumerate}
\end{athm}
\begin{aproof}
\argument{}{%
\ that for all $k\in\N\cap[(a\cgamma)^{\nicefrac{1}{\lambda}},\infty)$ it holds that
\begin{equation}\llabel{exp:nonnegative}
    0
    \leq
    a\gamma_k
    =
    a\cgamma k^{-\lambda}
    \leq1.
\end{equation}}
\argument{\lref{exp:nonnegative};
the fact that for all $x\in\R$ it holds that $1-x\leq\exp(-x)$;
the monotonicity of $x\mapsto x^{-\lambda}$}{%
that for all $m,n\in\N$ with $(a\cgamma)^{\nicefrac1\lambda} \leq m\leq n$ it holds that
\begin{equation}
\begin{split}
    \prod_{k=m}^n(1-a\gamma_k)
    &\leq
    \exp\left(
      -a\smallsum_{k=m}^n\gamma_k
    \right)
    \\
    &=
    \exp\left(
      -a\cgamma
      \smallsum_{k=m}^nk^{-\lambda}
    \right)
    \\
    &\leq
    \exp\left(
      -a\cgamma
      \int_m^{n+1}x^{-\lambda}\,dx
    \right)
    \\
    &=
    \exp\left(
      -\frac{a\cgamma}{1-\lambda}
      \left[
        (n+1)^{1-\lambda}
        -
        m^{1-\lambda}
      \right]
    \right).
\end{split}
\end{equation}}
This implies \cref{exp:1}.
For the proof of \cref{exp:2}, throughout the remainder of
this proof let $M\in(0,\infty)$ satisfy
\begin{equation}\llabel{exp:second:M}
    M
    =
    \max\left\{
      (a\cgamma)^{\nicefrac{1}{\lambda}},
      \left(
        \frac{4\lambda}{a\cgamma}
      \right)^{\nicefrac{1}{1-\lambda}}
    \right\}
\end{equation}
and for every $m,l\in\N$ with $M\leq m\leq l$, let
$S_{m,l}\in[0,\infty)$ satisfy
\begin{equation}\llabel{exp:second:S-def}
    S_{m,l}
    =
    \smallsum_{k=m}^l
    \gamma_k^2
    \prod_{j=k+1}^l
    (1-a\gamma_j).
\end{equation}
\argument{\lref{exp:nonnegative};
\lref{exp:second:M};
\lref{exp:second:S-def}}{%
that for all $m\in\N$ with $M\leq m$ it holds that
\begin{equation}\llabel{exp:second:base}
\begin{split}
    S_{m,m}
    &=
    \gamma_m^2
   \leq
    \frac1a\gamma_m
    \leq
    \frac2a\gamma_m.
\end{split}
\end{equation}}
\argument{the fact that for all $x\in(0,\infty)$ it holds that $(1+x)^\lambda\leq 1 +\lambda x$;
the fact that for all $l\in\N\cap[2,\infty)$ it holds that
$\frac{l}{l-1}\leq2$;
\lref{exp:second:M}}{%
that for all $m,l\in\N$ with $M\leq m<l$ it holds that
\begin{equation}\llabel{exp:second:step-difference}
\begin{split}
    \frac{\gamma_{l-1}-\gamma_l}{\gamma_l^2}
    &= \frac{l^{2\lambda}}{\cgamma^2}
    \left[
        \frac{\cgamma}{(l-1) ^\lambda}
      -
        \frac{\cgamma}{l ^\lambda}
    \right]
    \\
    &=
    \frac{l^\lambda}{\cgamma}
    \left[
      \left(
        \frac{l}{l-1}
      \right)^\lambda
      -1
    \right]
    \\
    &=
    \frac{l^\lambda}{\cgamma}
    \left[
      \left(
        1+\frac1{l-1}
      \right)^\lambda
      -1
    \right]
    \\
    &\leq
    \frac{\lambda l^\lambda}{\cgamma(l-1)}
    \\
    &=
    \frac{\lambda}{\cgamma}
    l^{\lambda-1}
    \frac{l}{l-1}
    \\
    &\leq
    \frac{2\lambda}{\cgamma}
    l^{\lambda-1}
    \\
    &=
    \frac{2\lambda}
    {\cgamma l^{1-\lambda}}
    \\
    &\leq
    \frac a2.
\end{split}
\end{equation}}
\argument{the fact that for all $l\in\N\cap[2,\infty)$ it holds
that $\gamma_{l-1}\geq\gamma_l$;
\lref{exp:second:step-difference}}{%
that for all $m,l\in\N$ with $M\leq m<l$ it holds that
\begin{equation}\llabel{exp:second:supersolution}
\begin{split}
    &\frac2a\gamma_l
    -
    \left[
      (1-a\gamma_l)
      \frac2a\gamma_{l-1}
      +
      \gamma_l^2
    \right]
    \\
    &=
    \frac2a
    \left[
      \gamma_l-\gamma_{l-1}
      +
      a\gamma_l\gamma_{l-1}
      -
      \frac a2\gamma_l^2
    \right]
    \\
    &\geq
    \frac2a
    \left[
      -\frac a2\gamma_l^2
      +
      a\gamma_l^2
      -
      \frac a2\gamma_l^2
    \right]
    \\
    &=0,
\end{split}
\end{equation}
and hence
\begin{equation}\llabel{exp:second:supersolution-conclusion}
    (1-a\gamma_l)
    \frac2a\gamma_{l-1}
    +
    \gamma_l^2
    \leq
    \frac2a\gamma_l.
\end{equation}}
\argument{\lref{exp:second:S-def}}{%
that for all $m,l\in\N$ with $M\leq m<l$ it holds that
\begin{equation}\llabel{exp:second:S-recursion}
\begin{split}
    S_{m,l}
    &=
    \sum_{k=m}^{l-1}
    \gamma_k^2
    \prod_{j=k+1}^l
    (1-a\gamma_j)
    +
    \gamma_l^2
    \\
    &=
    (1-a\gamma_l)
    \sum_{k=m}^{l-1}
    \gamma_k^2
    \prod_{j=k+1}^{l-1}
    (1-a\gamma_j)
    +
    \gamma_l^2
    \\
    &=
    (1-a\gamma_l)S_{m,l-1}
    +
    \gamma_l^2.
\end{split}
\end{equation}}
\argument{\lref{exp:nonnegative};
\lref{exp:second:base};
\lref{exp:second:supersolution-conclusion};
\lref{exp:second:S-recursion};
induction on $l-m$}{%
that for all $m,l\in\N$ with $M\leq m\leq l$ it holds that
\begin{equation}\llabel{exp:second:S-bound}
    S_{m,l}
    \leq
    \frac2a\gamma_l.
\end{equation}}
\argument{\lref{exp:second:M};
\lref{exp:second:S-def};
\lref{exp:second:S-bound};
the fact that for all $n\in\N$ it holds that
$\gamma_n=\cgamma n^{-\lambda}$}{%
that for all $m,n\in\N$ with
$\max\left\{
      (a\cgamma)^{\nicefrac{1}{\lambda}},
      \left(
        \frac{4\lambda}{a\cgamma}
      \right)^{\nicefrac{1}{1-\lambda}}
    \right\}
    \leq m\leq n$
it holds that
\begin{equation}
\begin{split}
    \sum_{k=m}^n
    \gamma_k^2
    \prod_{j=k+1}^n
    (1-a\gamma_j)
    &=
    S_{m,n}
    \leq
    \frac2a\gamma_n
    =
    \frac{2\cgamma}{a}
    n^{-\lambda}.
\end{split}
\end{equation}}
\end{aproof}

\begin{athm}{theorem}{theo:Vn}
    Assume \cref{setting2}, assume $\const < \infty$, let $c \in(0, \infty)$, $\vartheta \in \R^\fd$ satisfy for all $\theta \in \R^\fd$ that
    \begin{equation}
        \langle\theta-\vartheta, \E[\Gbatch(\theta, \Inn_{1,1})]\rangle \leq -c \max\{\|\theta-\vartheta\|^2,\|\E[\Gbatch(\theta,\Inn_{1,1})]\|^2\},
    \end{equation}
    let $\lambda\in (0,1)$, $\cgamma\in (0,\infty)$,
    assume for all $n\in\N$ that $\gamma_n=\cgamma n^{-\lambda}$,
    let $\widehat N \in \N$  satisfy $\widehat N=\left\lceil\max \left\{ 2, \big(\frac{2(1+\alpha)c\cgamma}{1-\alpha}\big)^{\nicefrac{1}{\lambda}}, \big(\frac{4\lambda}{(1-\alpha)c\cgamma}\big)^{\nicefrac{1}{1-\lambda}}, \big(\frac{2(2-\alpha)(J-1+\alpha)\cgamma}{(1-\alpha)^2c^3 J}\big)^{\nicefrac{1}{\lambda}} \right\}\right\rceil$,
    and for every $n \in \N_0$ let $V_n \in [0,\infty)$ satisfy \begin{equation}
        V_n = \E[\|\Theta_n - \vartheta\|^2] + \frac{1+\alpha}{1-\alpha} \gamma_n \max\{\E[\langle \Theta_{n} - \vartheta, \m_{n}\rangle],0\} + \frac{2}{(1-\alpha)^2} \gamma_n^2 \E[\|\m_n\|^2].
    \end{equation}
    Then 
    for all $n \in \N$ with $n \geq \widehat N$ it holds that
        \begin{equation}
        \begin{split}
             \E[\|\Theta_n - \vartheta\|^2] &\leq V_{\widehat N-1} \exp \left(\frac{(1-\alpha)c\cgamma}{1-\lambda}\widehat N^{1-\lambda}\right) \exp \left(-\frac{(1-\alpha)c\cgamma}{1-\lambda} n^{1-\lambda}\right) \\
             & \quad + \frac{4(2-\alpha)\cgamma \const^2}{(1-\alpha)c J} n^{-\lambda}.
        \end{split}
        \end{equation} 
\end{athm}
\begin{aproof}
Throughout this proof let $a \in (0,\infty)$, $b_J\in [0,\infty)$ satisfy
\begin{equation}
    a=(1-\alpha)c \qandq b_J=\frac{2(2-\alpha)\const^2}{J}.
\end{equation}
\argument{the definition of $\widehat N$;
the fact that for all $n\in\N$ it holds that $\gamma_n=\cgamma n^{-\lambda}$}{%
that for all $n\in\N\cap[\widehat N,\infty)$ it holds that
\begin{equation}\llabel{theo:Vn:small-one}
    \gamma_n
    \leq
    \frac{1-\alpha}{2(1+\alpha)c},
\end{equation}
\begin{equation}\llabel{theo:Vn:small-two}
    \gamma_n
    \leq
    \frac{
      (1-\alpha)^2Jc^3
    }{
      2(2-\alpha)(J-1+\alpha)
    },
\end{equation}
and
\begin{equation}\llabel{theo:Vn:large-index}
    \widehat N
    \geq
    \left(
      \frac{4\lambda}{a\cgamma}
    \right)^{1/(1-\lambda)}.
\end{equation}}
\argument{the fact that $(\gamma_n)_{n\in\N}$ is nonincreasing;
\lref{theo:Vn:small-one};
\cref{prop:Vn:recursionN}}{%
that for all $n\in\N\cap[\widehat N,\infty)$ it holds that
\begin{equation}\llabel{theo:Vn:beta-recursion}
    V_n
    \leq
    \beta_nV_{n-1}
    +
    b_J\gamma_n^2.
\end{equation}}
\argument{the definition of $\beta_n$;
\lref{theo:Vn:small-two}}{%
that for all $n\in\N\cap[\widehat N,\infty)$ it holds that
\begin{equation}\llabel{theo:Vn:beta-upper}
\begin{split}
    \beta_n
    &=
    1-2a\gamma_n
    +
    \frac{
      2(2-\alpha)(J-1+\alpha)
    }{
      (1-\alpha)Jc^2
    }
    \gamma_n^2
    \\
    &\leq
    1-a\gamma_n.
\end{split}
\end{equation}}
\argument{\lref{theo:Vn:beta-recursion};
\lref{theo:Vn:beta-upper};
the fact that for all $n\in \N$ it holds that $V_{n-1}\geq0$}{%
that for all $n\in\N\cap[\widehat N,\infty)$ it holds that
\begin{equation}\llabel{theo:Vn:scalar-recursion}
    V_n
    \leq
    (1-a\gamma_n)V_{n-1}
    +
    b_J\gamma_n^2.
\end{equation}}
\argument{\lref{theo:Vn:scalar-recursion};}{%
that for all $n\in\N\cap[\widehat N,\infty)$ it holds that
\begin{equation}\llabel{theo:Vn:unrolled}
\begin{split}
    V_n
    &\leq
    V_{\widehat N-1}
    \prod_{k=\widehat N}^n
    (1-a\gamma_k)
    \\
    &\quad
    +
    b_J
    \sum_{j=\widehat N}^n
    \gamma_j^2
    \prod_{k=j+1}^n
    (1-a\gamma_k).
\end{split}
\end{equation}}
\argument{\cref{exp};
\lref{theo:Vn:large-index};
\lref{theo:Vn:unrolled}}{%
that for all $n\in\N\cap[\widehat N,\infty)$ it holds that
\begin{equation}\llabel{theo:Vn:explicit}
\begin{split}
    V_n
    &\leq
    V_{\widehat N-1}
    \exp\left(
      -\frac{a\cgamma}{1-\lambda}
      \left[
        (n+1)^{1-\lambda}
        -
        \widehat N^{1-\lambda}
      \right]
    \right)
    \\
    &\quad
    +
    \frac{2b_J\cgamma}{a}
    n^{-\lambda}.
\end{split}
\end{equation}}
\argument{the definitions of $a$ and $b_J$;
the fact that for all $n\in \N$ it holds that
$\E[\|\Theta_n-\vartheta\|^2]\leq V_n$;
\lref{theo:Vn:explicit}}{%
that for all $n\in\N\cap[\widehat N,\infty)$ it holds that
\begin{equation}
\begin{split}
    \E[\|\Theta_n-\vartheta\|^2]
    &\leq
    V_{\widehat N-1}
    \exp\left(
      -\frac{(1-\alpha)c\cgamma}{1-\lambda}
      \left[
        (n+1)^{1-\lambda}
        -
        \widehat N^{1-\lambda}
      \right]
    \right)
    \\
    &\quad
    +
    \frac{
      4(2-\alpha)\cgamma\const^2
    }{
      (1-\alpha)cJ
    }
    n^{-\lambda}.
\end{split}
\end{equation}}
\argument{the fact that for all $n\in\N$ it holds that
$(n+1)^{1-\lambda}\geq n^{1-\lambda}$}{%
that for all $n\in\N\cap[\widehat N,\infty)$ it holds that
\begin{equation}
\begin{split}
    \E[\|\Theta_n-\vartheta\|^2]
    &\leq
    V_{\widehat N-1}
    \exp\left(
      \frac{(1-\alpha)c\cgamma}{1-\lambda}
      \widehat N^{1-\lambda}
    \right)
    \exp\left(
      -\frac{(1-\alpha)c\cgamma}{1-\lambda}
      n^{1-\lambda}
    \right)
    \\
    &\quad +
    \frac{
      4(2-\alpha)\cgamma\const^2
    }{
      (1-\alpha)cJ
    }
    n^{-\lambda}.
\end{split}
\end{equation}}
\end{aproof}

\begin{athm}{cor}{cor:main-result}
Let $(\Omega,\mathcal F,\Pp)$ be a probability space, let
$\fd,\dimm\in\N$, let
$\Inn_{n,j}\colon\Omega\to\R^\dimm$,
$(n,j)\in\N_0^2$, be \iid\ random variables, let
$\widehat\Gbatch\colon\R^\fd\times\R^\dimm\to\R^\fd$ be measurable,
let $\xi,\vartheta\in\R^\fd$, $\alpha\in[0,1)$,
$c,\cgamma, \const \in(0,\infty)$, $\lambda\in(0,1)$,
for every $J\in\N$ let
$\Theta^J\colon\N_0\times\Omega\to\R^\fd$ and
$\widehat\m^J\colon\N_0\times\Omega\to\R^\fd$ satisfy for all
$n\in\N$ that
\begin{equation}\label{cor:main-result;1}
    \Theta_0^J=\xi,
    \qquad
    \widehat\m_0^J=0,
    \qquad
    \Theta_n^J
    =
    \Theta_{n-1}^J
    -
    \cgamma n^{-\lambda}\widehat\m_n^J,
\end{equation}
\begin{equation}\label{cor:main-result;2}
    \andq \widehat\m_n^J
    =
    \alpha\widehat\m_{n-1}^J
    +
    (1-\alpha)
    \left[
        \frac1J
        \smallsum_{j=1}^J
        \widehat\Gbatch(\Theta_{n-1}^J,\Inn_{n,j})
    \right],
\end{equation}
 and
assume for all $\theta\in\R^\fd$ that 
$\sup_{J,n\in\N}
    \big(
        \E\bigl[
            \|\widehat\Gbatch(\Theta_{n-1}^J,\Inn_{0,1})\|^2
        \bigr]
    \big)^{\nicefrac12}
    \leq \const$,
    $\E[\|\widehat\Gbatch(\theta,\Inn_{0,1})\|]<\infty$,
    and \begin{equation}\llabel{cor:main-result:coercivity}
    \langle
        \theta-\vartheta,
        \E[\widehat\Gbatch(\theta,\Inn_{1,1})]
    \rangle
    \geq
    c\max\{
        \|\theta-\vartheta\|^2,
        \|\E[\widehat\Gbatch(\theta,\Inn_{1,1})]\|^2
    \}.
\end{equation}

Then there exists $\C\in(0,\infty)$ such that for all
$J,n\in\N$ it holds that
\begin{equation}
    \E[\|\Theta_n^J-\vartheta\|^2]
    \leq \C \bigl(\exp (-\C^{-1} n^{1-\lambda}) + J^{-1} n^{-\lambda} \bigr).
\end{equation}
\end{athm}
\begin{aproof}
Throughout this proof let $p,q,\eta,D\in(0,\infty)$ satisfy
\begin{equation}\llabel{cor:main-result:constants}
    p=\frac{1+\alpha}{1-\alpha},
    \qquad
    q=\frac{2}{(1-\alpha)^2},
    \qquad
    \eta
    =
    \frac{(1-\alpha)c\cgamma}{1-\lambda},
    \qquad
    D
    =
    \frac{
        4(2-\alpha)\cgamma\const^2
    }{
        (1-\alpha)c
    },
\end{equation}
let $\Gbatch\colon
\R^\fd\times\R^\dimm\to\R^\fd$  satisfy that
    $\widehat{\Gbatch}
    =
    -\Gbatch$,
for every $J\in\N$ let
$\m^J\colon\N_0\times\Omega\to\R^\fd$ satisfy that
$
    \widehat{\m}^J
    =
    -\m^J$,
let $(\gamma_n)_{n \in \N_0}\subseteq (0,\infty)$ satisfy for all $n \in \N$ that $\gamma_0=\cgamma$ and $\gamma_n=\cgamma n^{-\lambda}$,
for every $J\in\N$, $n\in\N_0$ let
$V_n^J\in[0,\infty)$ satisfy
\begin{equation}\llabel{cor:main-result:V}
\begin{split}
    V_n^J = \E[\|\Theta_n^J-\vartheta\|^2]+
    p\gamma_n
    \max\left\{
        \E[
            \langle
                \Theta_n^J-\vartheta,
                \m_n^J
            \rangle
        ],
        0
    \right\}
    +
    q\gamma_n^2
    \E[\|\m_n^J\|^2],
\end{split}
\end{equation}
let $\overline N\in\N$ satisfy
$
    \overline N
    =
    \Big\lceil
    \max\Big\{
        2,
        \left(
            \frac{
                2(1+\alpha)c\cgamma
            }{
                1-\alpha
            }
        \right)^{\nicefrac1\lambda},
        \left(
            \frac{
                4\lambda
            }{
                (1-\alpha)c\cgamma
            }
        \right)^{\nicefrac1{1-\lambda}},
        \left(
            \frac{
                2(2-\alpha)\cgamma
            }{
                (1-\alpha)^2c^3
            }
        \right)^{\nicefrac1\lambda}
    \Big\}
    \Big\rceil
$,
for every $J \in \N$ let $N_J \in \N$ satisfy $N_J =\Big\lceil
    \max\Big\{ 2,
        \left(
            \frac{
                2(1+\alpha)c\cgamma
            }{
                1-\alpha
            }
        \right)^{\nicefrac1\lambda},
        \left(
            \frac{
                4\lambda
            }{
                (1-\alpha)c\cgamma
            }
        \right)^{\nicefrac1{1-\lambda}},
        \left(
            \frac{
                2(2-\alpha)(J-1+\alpha)\cgamma
            }{
                (1-\alpha)^2c^3J
            }
        \right)^{\nicefrac1\lambda}
    \Big\}
    \Big\rceil$,
    and let $R\in[0,\infty)$ satisfy
\begin{equation}\llabel{cor:main-result:R}
    R
    =
    \|\xi-\vartheta\|
    +
    \cgamma\const
    \smallsum_{k=1}^{\overline N}k^{-\lambda}.
\end{equation}
\argument{\cref{cor:main-result;1};\cref{cor:main-result;2}}
{%
that for all $J,n\in\N$ it holds that
\begin{equation}\llabel{cor:main-result:transformed-recursion}
\begin{split}
    \Theta_n^J
    &=
    \Theta_{n-1}^J
    +
    \gamma_n\m_n^J,
    \qandq
    \m_n^J
    =
    \alpha\m_{n-1}^J
    +
    (1-\alpha)
    \left[
        \frac1J
        \smallsum_{j=1}^J
        \Gbatch(\Theta_{n-1}^J,\Inn_{n,j})
    \right].
\end{split}
\end{equation}}
\argument{\lref{cor:main-result:coercivity}}{%
that for all $\theta\in\R^\fd$ it holds that
\begin{equation}\llabel{cor:main-result:negative-coercivity}
\begin{split}
    \left\langle
        \theta-\vartheta,
        \E[
            \Gbatch(\theta,\Inn_{1,1})
        ]
    \right\rangle
    &=
    -\langle
        \theta-\vartheta,
        \E[
            \widehat{\Gbatch}(\theta,\Inn_{1,1})
        ]
    \rangle
    \\
    &\leq -c\max\{
        \|\theta-\vartheta\|^2,
        \|
            \E[
                \widehat{\Gbatch}(\theta,\Inn_{1,1})
            ]
        \|^2
    \}
    \\
    &\leq -c\max\{
        \|\theta-\vartheta\|^2,
        \|
            \E[
                \Gbatch(\theta,\Inn_{1,1})
            ]
        \|^2
    \}.
\end{split}
\end{equation}}
\argument{the fact that for all $J\in\N$ it holds that
$\frac{J-1+\alpha}{J}\leq1$}{%
that for all $J\in\N$ it holds that
\begin{equation}\llabel{cor:main-result:N-order}
\begin{split}
    N_J
    &\leq 
    \overline N.
    \end{split}
\end{equation}}
\argument{the momentum recursion;
the Minkowski inequality;
the fact that
$(\Inn_{n,j})_{(n,j)\in\N_0^2}$ are \iid;
the definition of $\const$}{%
that for all $J\in\N$, $n\in\N_0$ it holds that
\begin{equation}\llabel{cor:main-result:m-bound}
\begin{split}
    \left(
        \E[\|\m_n^J\|^2]
    \right)^{\nicefrac12}
    &\leq
    (1-\alpha)
    \smallsum_{k=1}^n
    \alpha^{n-k}
    \frac1J
    \smallsum_{j=1}^J
    \left(
        \E[
            \|\Gbatch(\Theta_{k-1}^J,\Inn_{k,j})\|^2
        ]
    \right)^{\nicefrac12}
    \\
    &\leq
    (1-\alpha)
    \smallsum_{k=1}^n
    \alpha^{n-k}\const
    \\
    &=
    (1-\alpha^n)\const
    \leq
    \const.
\end{split}
\end{equation}}
\argument{\cref{cor:main-result;1};
the Minkowski inequality;
\lref{cor:main-result:m-bound};
\lref{cor:main-result:R}}{%
that for all $J\in\N$,
$n\in\N_0\cap[0,\overline N]$ it holds that
\begin{equation}\llabel{cor:main-result:theta-bound}
\begin{split}
    \left(
        \E[\|\Theta_n^J-\vartheta\|^2]
    \right)^{\nicefrac12}
    &\leq
    \|\xi-\vartheta\|
    +
    \smallsum_{k=1}^n
    \gamma_k
    \left(
        \E[\|\m_k^J\|^2]
    \right)^{\nicefrac12}
    \\
    &\leq R.
\end{split}
\end{equation}}
\argument{the Cauchy--Schwarz inequality;
\lref{cor:main-result:V};
\lref{cor:main-result:m-bound};
\lref{cor:main-result:theta-bound}}{%
that for all $J\in\N$,
$n\in\N_0\cap[0,\overline N]$ it holds that
\begin{equation}\llabel{cor:main-result:V-bound}
    V_n^J\leq R^2
    +
    p\cgamma R\const
    +
    q\cgamma^2\const^2.
\end{equation}}
\argument{\lref{cor:main-result:transformed-recursion};
\lref{cor:main-result:negative-coercivity};
\lref{cor:main-result:N-order};
\lref{cor:main-result:V-bound};
\cref{theo:Vn}}{%
that for all $J,n\in\N$ with $n\geq \overline N$ it holds that
\begin{equation}\llabel{cor:main-result:large-n}
\begin{split}
    \E[\|\Theta_n^J-\vartheta\|^2]
    &\leq
    (R^2
    +
    p\cgamma R\const
    +
    q\cgamma^2\const^2)
    \exp\left(
        \eta\overline N^{1-\lambda}
    \right)
    \exp\left(
        -\eta n^{1-\lambda}
    \right)
    +
    \frac{D}{J}n^{-\lambda}.
\end{split}
\end{equation}}
Let $\C\in(0,\infty)$ satisfy
\begin{equation}\llabel{cor:main-result:C}
    \C
    =
    \max\left\{
        1,\,
        \eta^{-1},\,
        (R^2
    +
    p\cgamma R\const
    +
    q\cgamma^2\const^2)
        \exp\left(
            \eta\overline N^{1-\lambda}
        \right),\,
        D,\,
        R^2\exp\left(
            \overline N^{1-\lambda}
        \right)
    \right\}.
\end{equation}
\argument{\lref{cor:main-result:N-order};
\lref{cor:main-result:theta-bound};
\lref{cor:main-result:large-n};
\lref{cor:main-result:C}}{%
that for all $J,n\in\N$ it holds that
\begin{equation}
    \E[\|\Theta_n^J-\vartheta\|^2]
    \leq \C \bigl(\exp (-\C^{-1} n^{1-\lambda}) + J^{-1} n^{-\lambda} \bigr).
\end{equation}}
\end{aproof}

\subsection{Inverse-linear step sizes}\label{section:rates2}
We analyze the rate of convergence of the \SGD\ with momentum algorithm to the zero of the function $\f$ for inverse-linear step size sequences. Similarly to \cref{section:rates} we first calculate in the elementary result in \cref{exp:lambda-one} below convergence rates for each unrolled coefficient of our recursive bound and then establish the convergence rate in \cref{theo:Vnlambda} where we obtain three regimes.
Only for completeness we include here a proof for \cref{exp:lambda-one}.

\begin{athm}{lemma}{exp:lambda-one}
    Let $a,\cgamma\in(0,\infty)$, let
    $(\gamma_n)_{n\in\N}\subseteq(0,\infty)$ satisfy for all
    $n\in\N$ that $\gamma_n=\cgamma n^{-1}$,
    and let $\rho\in(0,\infty)$ satisfy $\rho=a\cgamma.$
    Then
    \begin{enumerate}[label=(\roman*)]
        \item \label{exp:lambda-one:1} for all $m,n\in\N$ with
        $\rho\leq m\leq n$ it holds that
        \begin{equation}
            \prod_{k=m}^n
            (1-a\gamma_k)
            \leq
            \left(
                \frac{m}{n+1}
            \right)^\rho,
        \end{equation}
        \item \label{exp:lambda-one:2} for all
        $m,n\in\N$ with $\max\{2,\rho\} \leq m\leq n$ it holds that
        \begin{equation}
        \sum_{k=m}^n
            \gamma_k^2
            \prod_{j=k+1}^n
            (1-a\gamma_j)
            \leq
            \begin{cases}
           \frac{
                2^\rho\cgamma^2
            }{
                (1-\rho)(m-1)^{1-\rho}
            }
            (n+1)^{-\rho} & \colon \rho<1, \\
             \frac{\cgamma^2}{n}
            \left[
                \frac1m+\log\left(\frac nm\right)
            \right] & \colon \rho = 1, \\
            \frac{\cgamma^2}{(\rho-1)n} & \colon \rho > 1.
    \end{cases}
        \end{equation}
    \end{enumerate}
\end{athm}
\begin{aproof}
\argument{the definitions of $\rho$ and $\gamma_k$}{%
that for all $k\in\N$ with $\rho\leq k$ it holds that
\begin{equation}\llabel{exp:lambda-one:nonnegative}
\begin{split}
    0
    &\leq
    a\gamma_k
    =
    \frac{a\cgamma}{k}
    =
    \frac{\rho}{k}
    \leq1.
\end{split}
\end{equation}}
\argument{\lref{exp:lambda-one:nonnegative};
the fact that for all $x\in\R$ it holds that
$1-x\leq\exp(-x)$;
the monotonicity of $x\mapsto x^{-1}$}{%
that for all $m,n\in\N$ with $\rho\leq m\leq n$ it holds that
\begin{equation}\llabel{exp:lambda-one:product}
\begin{split}
    \prod_{k=m}^n(1-a\gamma_k)
    &\leq
    \exp\left(
        -a\sum_{k=m}^n\gamma_k
    \right)
    \\
    &=
    \exp\left(
        -\rho\sum_{k=m}^n\frac1k
    \right)
    \\
    &\leq
    \exp\left(
        -\rho\int_m^{n+1}\frac1x\,dx
    \right)
    \\
    &=
    \exp\left(
        -\rho
        \log\left(\frac{n+1}{m}\right)
    \right)
    \\
    &=
    \left(
        \frac{m}{n+1}
    \right)^\rho.
\end{split}
\end{equation}}
This proves \cref{exp:lambda-one:1}.
For the proof of \cref{exp:lambda-one:2} we distinguish between the case $\rho<1$, the case $\rho=1$, and the case $\rho>1$. We first show \cref{exp:lambda-one:2} in the case
    \begin{equation}\llabel{exp:hp1}
        \rho<1.
    \end{equation}
    \startnewargseq
\argument{\lref{exp:lambda-one:product};
}{%
that for all $m,n,k\in\N$ with
$2\leq m\leq k\leq n$ it holds that
\begin{equation}\llabel{exp:lambda-one:subproduct}
    \prod_{j=k+1}^n
    (1-a\gamma_j)
    \leq
    \left(
        \frac{k+1}{n+1}
    \right)^\rho.
\end{equation}}
\argument{the fact that for all $k\in\N$ it holds that
$k+1\leq2k$;
\lref{exp:lambda-one:subproduct}}{%
that for all $m,n\in\N$ with $2\leq m\leq n$ it holds that
\begin{equation}\llabel{exp:lambda-one:subcritical-first}
\begin{split}
    &\sum_{k=m}^n
    \gamma_k^2
    \prod_{j=k+1}^n
    (1-a\gamma_j)
    \leq
    \frac{\cgamma^2}{(n+1)^\rho}
    \sum_{k=m}^n
    \frac{(k+1)^\rho}{k^2}
    \leq
    \frac{2^\rho\cgamma^2}{(n+1)^\rho}
    \sum_{k=m}^n
    k^{\rho-2}.
\end{split}
\end{equation}}
\argument{the fact that the function
$x\mapsto x^{\rho-2}$ is decreasing;
\lref{exp:hp1}}{%
that for all $m,n\in\N$ with $2\leq m\leq n$ it holds that
\begin{equation}\llabel{exp:lambda-one:subcritical-sum}
\begin{split}
    \sum_{k=m}^n
    k^{\rho-2}
    &\leq
    \int_{m-1}^n
    x^{\rho-2}\,dx
    \\
    &=
    \frac{
        (m-1)^{\rho-1}
        -
        n^{\rho-1}
    }{
        1-\rho
    }
    \\
    &\leq
    \frac{
        (m-1)^{\rho-1}
    }{
        1-\rho
    }
    \\
    &=
    \frac1{
        (1-\rho)(m-1)^{1-\rho}
    }.
\end{split}
\end{equation}}
\argument{\lref{exp:lambda-one:subcritical-first};
\lref{exp:lambda-one:subcritical-sum}}{%
that for all $m,n\in\N$ with $2\leq m\leq n$ it holds that
\begin{equation}
\begin{split}
    &\sum_{k=m}^n
    \gamma_k^2
    \prod_{j=k+1}^n
    (1-a\gamma_j)
    \\
    &\leq
    \frac{
        2^\rho\cgamma^2
    }{
        (1-\rho)(m-1)^{1-\rho}
    }
    (n+1)^{-\rho}.
\end{split}
\end{equation}}
We now consider the case
    \begin{equation}\llabel{exp:hp2}
        \rho=1.
    \end{equation}
    \startnewargseq
 \startnewargseq
\argument{the fact that for all $j\in\N$ it holds that
$a\gamma_j=\frac1j$}{%
that for all $m,n,k\in\N$ with $m\leq k\leq n$ it holds that
\begin{equation}\llabel{exp:lambda-one:critical-product}
\begin{split}
    \prod_{j=k+1}^n
    (1-a\gamma_j)
    &=
    \prod_{j=k+1}^n
    \left(
        1-\frac1j
    \right)
   =
    \prod_{j=k+1}^n
    \frac{j-1}{j}
    =
    \frac{k}{n}.
\end{split}
\end{equation}}
\argument{\lref{exp:lambda-one:critical-product};
the fact that for all $k\in\N$ it holds that $\gamma_k=\cgamma k^{-1}$}{%
that for all $m,n\in\N$ with $m\leq n$ it holds that
\begin{equation}\llabel{exp:lambda-one:critical-first}
\begin{split}
    &\sum_{k=m}^n
    \gamma_k^2
    \prod_{j=k+1}^n
    (1-a\gamma_j)
    =
    \sum_{k=m}^n
    \frac{\cgamma^2}{k^2}
    \frac{k}{n}
    =
    \frac{\cgamma^2}{n}
    \sum_{k=m}^n\frac1k.
\end{split}
\end{equation}}
\argument{the monotonicity of the function $x\mapsto x^{-1}$}{%
that for all $m,n\in\N$ with $m\leq n$ it holds that
\begin{equation}\llabel{exp:lambda-one:harmonic}
\begin{split}
    \sum_{k=m}^n\frac1k
    &\leq
    \frac1m
    +
    \int_m^n\frac1x\,dx
    =
    \frac1m
    +
    \log\left(\frac nm\right).
\end{split}
\end{equation}}
\argument{\lref{exp:lambda-one:critical-first};
\lref{exp:lambda-one:harmonic}}{%
that for all $m,n\in\N$ with $m\leq n$ it holds that
\begin{equation}
\begin{split}
    &\sum_{k=m}^n
    \gamma_k^2
    \prod_{j=k+1}^n
    (1-a\gamma_j)
    \leq
    \frac{\cgamma^2}{n}
    \left[
        \frac1m
        +
        \log\left(\frac nm\right)
    \right].
\end{split}
\end{equation}}
Finally, we consider the case
    \begin{equation}\llabel{exp:hp3}
        \rho>1.
    \end{equation}
\startnewargseq
For every $m,l\in\N$ with $\rho\leq m\leq l$, let
$S_{m,l}\in[0,\infty)$ satisfy
\begin{equation}\llabel{exp:lambda-one:S-def}
    S_{m,l}
    =
    \sum_{k=m}^l
    \gamma_k^2
    \prod_{j=k+1}^l
    (1-a\gamma_j).
\end{equation}
\argument{\lref{exp:lambda-one:S-def};
the fact that $\rho-1\leq m$}{%
that for all $m\in\N$ with $\rho\leq m$ it holds that
\begin{equation}\llabel{exp:lambda-one:supercritical-base}
\begin{split}
    S_{m,m}
    &=
    \gamma_m^2
    =
    \frac{\cgamma^2}{m^2}
    \leq
    \frac{\cgamma^2}{(\rho-1)m}.
\end{split}
\end{equation}}
\argument{}{%
for all $m,l\in\N$ with $\rho\leq m<l$ it holds that
\begin{equation}\llabel{exp:lambda-one:supercritical-super}
\begin{split}
    &\frac{\cgamma^2}{(\rho-1)l}
    -
    \left[
        \left(
            1-\frac{\rho}{l}
        \right)
        \frac{\cgamma^2}{(\rho-1)(l-1)}
        +
        \frac{\cgamma^2}{l^2}
    \right]
    =
    \frac{\cgamma^2}{l(l-1)}
    -
    \frac{\cgamma^2}{l^2}
    =
    \frac{\cgamma^2}{l^2(l-1)}
    \geq0.
\end{split}
\end{equation}}
\argument{\lref{exp:lambda-one:S-def}}{%
that for all $m,l\in\N$ with $\rho\leq m<l$ it holds that
\begin{equation}\llabel{exp:lambda-one:S-recursion}
\begin{split}
    S_{m,l}
    &=
    \left(
        1-a\gamma_l
    \right)
    S_{m,l-1}
    +
    \gamma_l^2
    =
    \left(
        1-\frac{\rho}{l}
    \right)
    S_{m,l-1}
    +
    \frac{\cgamma^2}{l^2}.
\end{split}
\end{equation}}
\argument{\lref{exp:lambda-one:nonnegative};
\lref{exp:lambda-one:supercritical-base};
\lref{exp:lambda-one:supercritical-super};
\lref{exp:lambda-one:S-recursion};
induction on $l-m$}{%
that for all $m,l\in\N$ with $\rho\leq m\leq l$ it holds that
\begin{equation}\llabel{exp:lambda-one:S-bound}
    S_{m,l}
    \leq
    \frac{\cgamma^2}{(\rho-1)l}.
\end{equation}}
\argument{\lref{exp:lambda-one:S-def};
\lref{exp:lambda-one:S-bound}}{%
that for all $m,n\in\N$ with $\rho\leq m\leq n$ it holds that
\begin{equation}
    \sum_{k=m}^n
    \gamma_k^2
    \prod_{j=k+1}^n
    (1-a\gamma_j)
    \leq
    \frac{\cgamma^2}{(\rho-1)n}.
\end{equation}}
\end{aproof}

\begin{athm}{theorem}{theo:Vnlambda}
    Assume \cref{setting2}, assume $\const < \infty$, let $c \in(0, \infty)$, $\vartheta \in \R^\fd$ satisfy for all $\theta \in \R^\fd$ that
    \begin{equation}
        \langle\theta-\vartheta, \E[\Gbatch(\theta, \Inn_{1,1})]\rangle \leq -c \max\{\|\theta-\vartheta\|^2,\|\E[\Gbatch(\theta,\Inn_{1,1})]\|^2\},
    \end{equation}
    let $\cgamma\in (0,\infty)$,
    assume for all $n\in\N$ that $\gamma_n=\cgamma n^{-1}$,
    let $\rho\in(0,\infty)$ satisfy $\rho=(1-\alpha)c\cgamma$,
    let $\widehat N \in \N$ satisfy $\widehat N= \left\lceil\max \left\{2,\big(\frac{2(1+\alpha)c\cgamma}{1-\alpha}\big), \big(\frac{2(2-\alpha)(J-1+\alpha)\cgamma}{(1-\alpha)^2c^3 J}\big) \right\}\right\rceil$,
    and for every $n \in \N_0$ let $V_n \in [0,\infty)$ satisfy \begin{equation}
        V_n = \E[\|\Theta_n - \vartheta\|^2] + \frac{1+\alpha}{1-\alpha} \gamma_n \max\{\E[\langle \Theta_{n} - \vartheta, \m_{n}\rangle],0\} + \frac{2}{(1-\alpha)^2} \gamma_n^2 \E[\|\m_n\|^2].
    \end{equation}
    Then for all
        $n\in\N\cap[2,\infty)$ with $n\geq \widehat N$ it holds that
        \begin{equation}
            \E[\|\Theta_n - \vartheta\|^2] \leq 
            \begin{cases}
            V_{\widehat N -1}
            \left(
                \frac{\widehat N}{n+1}
            \right)^\rho
           +
            \left(\frac{
                2^{\rho+1}\cgamma^2 (2-\alpha)\const^2
            }{
                (1-\rho)(\widehat N-1)^{1-\rho}
            }\right)
            \frac{1}{J (n+1)^{\rho}} & \colon \rho<1, \\
            V_{\widehat N-1}\left(
            \frac{\widehat N}{n+1}\right)
            +
            (2\cgamma^2(2-\alpha)\const^2) 
            \left[
                \frac{1}{\widehat N}
                +
                \log\left(\frac{ n}{\widehat N}\right)
            \right] \frac{1}{Jn}
            & \colon \rho = 1, \\
            V_{\widehat N-1}
            \left(
                \frac{\widehat N}{n+1}
            \right)^\rho
            + \left(\frac{2\cgamma^2(2-\alpha)\const^2}{(\rho-1)}\right)
            \frac{1}{Jn} & \colon \rho > 1.
    \end{cases}
        \end{equation}
\end{athm}
\begin{aproof}
Throughout this proof let $a\in(0,\infty)$ and
$b_J\in[0,\infty)$ satisfy
\begin{equation}\llabel{theo:Vnlambda:a-b}
    a=(1-\alpha)c
    \qquad\text{and}\qquad
    b_J=\frac{2(2-\alpha)\const^2}{J}.
\end{equation}
\argument{the definition of $\rho$;
\lref{theo:Vnlambda:a-b}}{%
that
\begin{equation}\llabel{theo:Vnlambda:rho}
    \rho
    =
    (1-\alpha)c\cgamma
    =
    a\cgamma.
\end{equation}}
\argument{the assumption that for all $n\in\N$ it holds that
$\gamma_n=\cgamma n^{-1}$}{%
that for all $n\in\N$ it holds that
\begin{equation}\llabel{theo:Vnlambda:monotone}
\begin{split}
    \gamma_{n+1}
    &=
    \frac{\cgamma}{n+1}
    \leq
    \frac{\cgamma}{n}
    =
    \gamma_n.
\end{split}
\end{equation}}
\argument{the definition of $\widehat N$;
}{%
that for all $n\in\N$ with $\widehat N\leq n$ it holds that
\begin{equation}\llabel{theo:Vnlambda:stepsize-small}
\begin{split}
    \gamma_n
    &=
    \frac{\cgamma}{n}
    \leq
    \frac{\cgamma}{\widehat N}
    \leq
    \frac{1-\alpha}{2(1+\alpha)c}.
\end{split}
\end{equation}}
\argument{\cref{prop:Vn:recursionN};
\lref{theo:Vnlambda:monotone};
\lref{theo:Vnlambda:stepsize-small}}{%
that for all $n\in\N$ with $\widehat N\leq n$ it holds that
\begin{equation}\llabel{theo:Vnlambda:beta-recursion}
    V_n
    \leq
    \beta_nV_{n-1}
    +
    b_J\gamma_n^2.
\end{equation}}
\argument{the definition of $\widehat N$}{%
that for all $n\in\N$ with $\widehat N\leq n$ it holds that
\begin{equation}\llabel{theo:Vnlambda:quadratic-bound}
\begin{split}
    &
    \frac{
      2(2-\alpha)(J-1+\alpha)
    }{
      (1-\alpha)Jc^2
    }
    \gamma_n
    \\
    &=
    \frac{
      2(2-\alpha)(J-1+\alpha)\cgamma
    }{
      (1-\alpha)Jc^2n
    }
    \\
    &\leq
    \frac{
      2(2-\alpha)(J-1+\alpha)\cgamma
    }{
      (1-\alpha)Jc^2\widehat N
    }
    \leq
    (1-\alpha)c
    \\
    &=a.
\end{split}
\end{equation}}
\argument{the definition of $\beta_n$;
\lref{theo:Vnlambda:a-b};
\lref{theo:Vnlambda:quadratic-bound}}{%
that for all $n\in\N$ with $\widehat N\leq n$ it holds that
\begin{equation}\llabel{theo:Vnlambda:beta-upper}
\begin{split}
    \beta_n
    &=
    1-2(1-\alpha)c\gamma_n
    +
    \frac{
      2(2-\alpha)(J-1+\alpha)
    }{
      (1-\alpha)Jc^2
    }
    \gamma_n^2
    \\
    &=
    1-2a\gamma_n
    +
    \frac{
      2(2-\alpha)(J-1+\alpha)
    }{
      (1-\alpha)Jc^2
    }
    \gamma_n^2
    \\
    &\leq
    1-2a\gamma_n+a\gamma_n
    \\
    &=
    1-a\gamma_n.
\end{split}
\end{equation}}
\argument{\lref{theo:Vnlambda:beta-recursion};
\lref{theo:Vnlambda:beta-upper};
the fact that for all $n\in\N_0$ it holds that $V_n\geq0$}{%
that for all $n\in\N$ with $\widehat N\leq n$ it holds that
\begin{equation}\llabel{theo:Vnlambda:scalar-recursion}
    V_n
    \leq
    (1-a\gamma_n)V_{n-1}
    +
    b_J\gamma_n^2.
\end{equation}}
\argument{the definition of $\widehat N$;
the fact that $\alpha\in[0,1)$;
\lref{theo:Vnlambda:rho}}{%
that
\begin{equation}\llabel{theo:Vnlambda:N-rho}
\begin{split}
    \widehat N
    &\geq
    \frac{2(1+\alpha)c\cgamma}{1-\alpha}
    =
    \frac{2(1+\alpha)}{(1-\alpha)^2}
    \rho
    \geq
    \rho.
\end{split}
\end{equation}
}
\argument{\lref{theo:Vnlambda:rho};
\lref{theo:Vnlambda:N-rho}}{%
that for all $n\in\N$ with $\widehat N\leq n$ it holds that
\begin{equation}\llabel{theo:Vnlambda:factor-nonnegative}
\begin{split}
    0
    &\leq
    a\gamma_n
    =
    \frac{a\cgamma}{n}
    =
    \frac{\rho}{n}
    \leq1,
\end{split}
\end{equation}
and hence
\begin{equation}
    0\leq1-a\gamma_n\leq1.
\end{equation}}
\argument{\lref{theo:Vnlambda:scalar-recursion};
\lref{theo:Vnlambda:factor-nonnegative};
induction on $n$}{%
that for all $n\in\N$ with $\widehat N\leq n$ it holds that
\begin{equation}\llabel{theo:Vnlambda:unrolled}
\begin{split}
    V_n
    &\leq
    V_{\widehat N-1}
    \prod_{k=\widehat N}^n
    (1-a\gamma_k)
    \\
    &\quad
    +
    b_J
    \smallsum_{k=\widehat N}^n
    \gamma_k^2
    \prod_{j=k+1}^n
    (1-a\gamma_j).
\end{split}
\end{equation}}
\argument{\cref{exp:lambda-one:1} in \cref{exp:lambda-one};
\lref{theo:Vnlambda:rho};
\lref{theo:Vnlambda:N-rho}}{%
that for all $n\in\N$ with $\widehat N\leq n$ it holds that
\begin{equation}\llabel{theo:Vnlambda:product-bound}
    \prod_{k=\widehat N}^n
    (1-a\gamma_k)
    \leq
    \left(
      \frac{\widehat N}{n+1}
    \right)^\rho.
\end{equation}}
\argument{\cref{exp:lambda-one:2} in \cref{exp:lambda-one};
\lref{theo:Vnlambda:rho};
\lref{theo:Vnlambda:N-rho}}{%
that for all $n\in\N$ with $\widehat N\leq n$ it holds that
\begin{equation}\llabel{theo:Vnlambda:convolution-bound}
\begin{split}
    &
    \sum_{k=\widehat N}^n
    \gamma_k^2
    \prod_{j=k+1}^n
    (1-a\gamma_j)
    \leq
    \begin{cases}
    \dfrac{
      2^\rho\cgamma^2
    }{
      (1-\rho)(\widehat N-1)^{1-\rho}
    }
    (n+1)^{-\rho}
    &\colon \rho<1,
    \\[3mm]
    \dfrac{\cgamma^2}{n}
    \left[
      \dfrac1{\widehat N}
      +
      \log\left(
        \dfrac n{\widehat N}
      \right)
    \right]
    &\colon \rho=1,
    \\[3mm]
    \dfrac{\cgamma^2}{(\rho-1)n}
    &\colon \rho>1.
    \end{cases}
\end{split}
\end{equation}}
\argument{\lref{theo:Vnlambda:unrolled};
\lref{theo:Vnlambda:product-bound};
\lref{theo:Vnlambda:convolution-bound}}{%
that for all $n\in\N$ with $\widehat N\leq n$ it holds that
\begin{equation}\llabel{theo:Vnlambda:V-bound}
\begin{split}
    V_n
    &\leq
    V_{\widehat N-1}
    \left(
      \frac{\widehat N}{n+1}
    \right)^\rho
   +
    b_J
    \begin{cases}
    \dfrac{
      2^\rho\cgamma^2
    }{
      (1-\rho)(\widehat N-1)^{1-\rho}
    }
    (n+1)^{-\rho}
    &\colon \rho<1,
    \\[3mm]
    \dfrac{\cgamma^2}{n}
    \left[
      \dfrac1{\widehat N}
      +
      \log\left(
        \dfrac n{\widehat N}
      \right)
    \right]
    &\colon \rho=1,
    \\[3mm]
    \dfrac{\cgamma^2}{(\rho-1)n}
    &\colon \rho>1.
    \end{cases}
\end{split}
\end{equation}}
\argument{the definition of $V_n$;
the fact that the second and third terms in the definition of
$V_n$ are nonnegative}{%
that for all $n\in\N$ it holds that
\begin{equation}\llabel{theo:Vnlambda:error-V}
    \E[\|\Theta_n-\vartheta\|^2]
    \leq
    V_n.
\end{equation}}
\argument{\lref{theo:Vnlambda:a-b};
\lref{theo:Vnlambda:V-bound};
\lref{theo:Vnlambda:error-V}}{%
that for all $n\in\N$ with $\widehat N\leq n$ it holds that
\begin{equation}
\begin{split}
    \E[\|\Theta_n-\vartheta\|^2]
    &\leq
    \begin{cases}
    V_{\widehat N-1}
    \left(
      \dfrac{\widehat N}{n+1}
    \right)^\rho
    +
    \left(
      \dfrac{
        2^{\rho+1}\cgamma^2
        (2-\alpha)\const^2
      }{
        (1-\rho)(\widehat N-1)^{1-\rho}
      }
    \right)
    \dfrac1{J(n+1)^\rho}
    &\colon \rho<1,
    \\[4mm]
    V_{\widehat N-1}
    \left(
      \dfrac{\widehat N}{n+1}
    \right)
    +
    \left(
      2\cgamma^2(2-\alpha)\const^2
    \right)
    \left[
      \dfrac1{\widehat N}
      +
      \log\left(
        \dfrac n{\widehat N}
      \right)
    \right]
    \dfrac1{Jn}
    &\colon \rho=1,
    \\[4mm]
    V_{\widehat N-1}
    \left(
      \dfrac{\widehat N}{n+1}
    \right)^\rho
    +
    \left(
      \dfrac{
        2\cgamma^2(2-\alpha)\const^2
      }{
        \rho-1
      }
    \right)
    \dfrac1{Jn}
    &\colon \rho>1.
    \end{cases}
\end{split}
\end{equation}}
\end{aproof}

\subsection{Convergence for general step sizes}\label{section:generalgamma}
\label{section:conv}
In this subsection we demonstrate the convergence of the \SGD\ with momentum algorithm to the zero of the function $\f$ for
general nonincreasing step size sequences which converge to zero and whose sum
diverges in \cref{conj:proof1}.

\begin{athm}{cor}{conj:proof1}
    Assume \cref{setting2}, let $c\in(0,\infty)$,
    $\vartheta\in\R^\fd$ satisfy for all $\theta\in\R^\fd$ that
    \begin{equation}\llabel{conj:proof1:coerc}
        \left\langle
        \theta-\vartheta,
        \E[\Gbatch(\theta,\Inn_{1,1})]
        \right\rangle
        \leq
        -c\max\left\{
        \|\theta-\vartheta\|^2,
        \|\E[\Gbatch(\theta,\Inn_{1,1})]\|^2
        \right\},
    \end{equation}
    assume $\const<\infty$, and assume for all $j\in\N$ that
    $\gamma_{j+1}\leq\gamma_j$, $\limsup\nolimits_{n\to\infty}\gamma_n=0$, and $\sum_{n=1}^\infty\gamma_n=\infty.$
    Then
    \begin{equation}
        \limsup\nolimits_{n\to\infty}
        \E[\|\Theta_n-\vartheta\|^2]
        =
        0.
    \end{equation}
\end{athm}
\begin{aproof}
Throughout this proof let $a \in(0,\infty)$ and
$b_J,d_J\in[0,\infty)$ satisfy
\begin{equation}\llabel{conj:proof1:constants}
    a=(1-\alpha)c,
    \qquad
    b_J = \frac{2(2-\alpha)\const^2}{J},
    \qandq
    d_J=
    \frac{
        2(2-\alpha)(J-1+\alpha)
    }{
        (1-\alpha)Jc^2
    },
\end{equation}
for every $n\in\N$, let $\beta_n\in\R$ satisfy
\begin{equation}\llabel{conj:proof1:beta}
    \beta_n
    =
    1-2a\gamma_n+d_J\gamma_n^2,
\end{equation}
and for every $n\in\N$ let $V_n\in[0,\infty)$ satisfy
\begin{equation}\llabel{conj:proof1:V-def}
\begin{split}
    V_n
    &=
    \E[\|\Theta_n-\vartheta\|^2]+
    \frac{1+\alpha}{1-\alpha}
    \gamma_n
    \max\left\{
        \E[
        \langle\Theta_n-\vartheta,\m_n\rangle
        ],
        0
    \right\}+
    \frac{2}{(1-\alpha)^2}
    \gamma_n^2
    \E[\|\m_n\|^2].
\end{split}
\end{equation}
\argument{\lref{conj:proof1:V-def};
the nonnegativity of the second and third terms in the definition
of $V_n$}{%
that for all $n\in\N$ it holds that
\begin{equation}\llabel{conj:proof1:error-V}
    0
    \leq
    \E[\|\Theta_n-\vartheta\|^2]
    \leq
    V_n.
\end{equation}}
\argument{the fact that
$\limsup_{n\to\infty}\gamma_n=0$;
\lref{conj:proof1:constants}}{%
that there exists $M\in\N\cap[2,\infty)$ such that for all
$n\in\N\cap[M+1,\infty)$ it holds that
\begin{equation}\llabel{conj:proof1:small-one}
    \gamma_n
    \leq
    \frac{1-\alpha}{2(1+\alpha)c},
\end{equation}
\begin{equation}\llabel{conj:proof1:small-two}
    d_J \gamma_n
    \leq
    a,
\end{equation}
and
\begin{equation}\llabel{conj:proof1:nonnegative}
    a\gamma_n
    \in[0,1).
\end{equation}}
\argument{
\lref{conj:proof1:small-one};
the assumption that for all $n\in\N$ it holds that
$\gamma_{n+1}\leq\gamma_n$;
\cref{prop:Vn:recursionN}}{%
that for all $n\in\N\cap[M+1,\infty)$ it holds that
\begin{equation}\llabel{conj:proof1:beta-recursion}
    V_n
    \leq
    \beta_nV_{n-1}
    + b_J \gamma_n^2.
\end{equation}}
\argument{\lref{conj:proof1:beta};
\lref{conj:proof1:small-two}}{%
that for all $n\in\N\cap[M+1,\infty)$ it holds that
\begin{equation}\llabel{conj:proof1:beta-upper}
\begin{split}
    \beta_n
    &=
    1-2a\gamma_n
    +d_J\gamma_n^2
    \\
    &\leq
    1-2a\gamma_n
    +a\gamma_n
    \\
    &=
    1-a\gamma_n.
\end{split}
\end{equation}}
\argument{\lref{conj:proof1:beta-recursion};
\lref{conj:proof1:beta-upper};
the nonnegativity of $V_{n-1}$}{%
that for all $n\in\N\cap[M+1,\infty)$ it holds that
\begin{equation}\llabel{conj:proof1:rec}
    V_n
    \leq
    (1-a\gamma_n)V_{n-1}
    +
    b_J\gamma_n^2.
\end{equation}}
\argument{\lref{conj:proof1:rec};
\lref{conj:proof1:nonnegative};
induction on $n$}{%
that for all $n\in\N\cap[M+1,\infty)$ it holds that
\begin{equation}\llabel{conj:proof1:unroll}
\begin{split}
    V_n
    &\leq
    \left(
        \prod_{i=M+1}^{n}
        (1-a\gamma_i)
    \right)
    V_M+
    b_J
    \sum_{j=M+1}^{n}
    \gamma_j^2
    \prod_{i=j+1}^{n}
    (1-a\gamma_i).
\end{split}
\end{equation}}
\argument{the fact that for all $x\in\R$ it holds that
$1+x\leq\exp(x)$;
the fact that
$\sum_{n=1}^\infty\gamma_n=\infty$}{%
that for all $M'\in\N\cap[M+1,\infty)$ it holds that
\begin{equation}\llabel{conj:proof1:prod}
\begin{split}
    \limsup\nolimits_{n\to\infty}
    \prod_{i=M'}^n
    (1-a\gamma_i)
    &\leq
    \limsup\nolimits_{n\to\infty}
    \exp\left(
        -a
        \smallsum_{i=M'}^n\gamma_i
    \right)=0.
\end{split}
\end{equation}}
\argument{the fact that for all $j\in\N$ and
$n\in\N\cap[j,\infty)$ it holds that
\begin{equation}
    a\gamma_j
    \prod_{i=j+1}^n
    (1-a\gamma_i)
    =
    \prod_{i=j+1}^n
    (1-a\gamma_i)
    -
    \prod_{i=j}^n
    (1-a\gamma_i)
\end{equation};
\lref{conj:proof1:nonnegative}}{%
that for all $M'\in\N\cap[M+1,\infty)$ and
$n\in\N\cap[M',\infty)$ it holds that
\begin{equation}\llabel{conj:proof1:tail}
\begin{split}
    &b_J
    \sum_{j=M'}^n
    \gamma_j^2
    \prod_{i=j+1}^n
    (1-a\gamma_i)
    \\
    &\leq
    \frac{b_J}{a}
    \left(
        \sup\nolimits_{j\geq M'}\gamma_j
    \right)
    \sum_{j=M'}^n
    \left[
        \prod_{i=j+1}^n
        (1-a\gamma_i)
        -
        \prod_{i=j}^n
        (1-a\gamma_i)
    \right]
    \\
    &=
    \frac{b_J}{a}
    \left(
        \sup\nolimits_{j\geq M'}\gamma_j
    \right)
    \left[
        \prod_{i=n+1}^n
        (1-a\gamma_i)
        -
        \prod_{i=M'}^n
        (1-a\gamma_i)
    \right]
    \\
    &=
    \frac{b_J}{a}
    \left(
        \sup\nolimits_{j\geq M'}\gamma_j
    \right)
    \left[
        1-
        \prod_{i=M'}^n
        (1-a\gamma_i)
    \right]
    \\
    &\leq
    \frac{b_J}{a}
    \sup\nolimits_{j\geq M'}\gamma_j.
\end{split}
\end{equation}}
\argument{the fact that for all
$M'\in\N\cap[M+1,\infty)$,
$j\in\N\cap[M,M'-1]$, and
$n\in\N\cap[M',\infty)$ it holds that
\begin{equation}
    \prod_{i=j+1}^n
    (1-a\gamma_i)
    \leq
    \prod_{i=M'}^n
    (1-a\gamma_i)
\end{equation};
\lref{conj:proof1:prod}}{%
that for all $M'\in\N\cap[M+1,\infty)$ it holds that
\begin{equation}\llabel{conj:proof1:head}
\begin{split}
    &\limsup\nolimits_{n\to\infty}
    b_J
    \sum_{j=M+1}^{M'-1} 
    \gamma_j^2
    \prod_{i=j+1}^n
    (1-a\gamma_i)
    \\
    &\leq
    b_J
    \left(
        \sum_{j=M+1}^{M'-1}
        \gamma_j^2
    \right)
    \limsup\nolimits_{n\to\infty}
    \prod_{i=M'}^n
    (1-a\gamma_i)
    \\
    &=0.
\end{split}
\end{equation}}
\argument{\lref{conj:proof1:unroll};
\lref{conj:proof1:prod};
\lref{conj:proof1:tail};
\lref{conj:proof1:head}}{%
that for all $M'\in\N\cap[M+1,\infty)$ it holds that
\begin{equation}\llabel{conj:proof1:limsup-V}
    \limsup\nolimits_{n\to\infty}V_n
    \leq
   \frac{b_J}{a}
    \sup\nolimits_{j\geq M'}\gamma_j.
\end{equation}}
\argument{\lref{conj:proof1:limsup-V};
the fact that
\begin{equation}
     \lim\nolimits_{M'\to\infty}
    \sup\nolimits_{j\geq M'}\gamma_j
    =
    \limsup\nolimits_{j\to\infty}\gamma_j
    =
    0
\end{equation};
the nonnegativity of $V_n$}{%
that
\begin{equation}\llabel{conj:proof1:V-limit}
    \limsup\nolimits_{n\to\infty}V_n
    =
    0.
\end{equation}}
\argument{\lref{conj:proof1:error-V};
\lref{conj:proof1:V-limit};
the nonnegativity of
$\E[\|\Theta_n-\vartheta\|^2]$}{%
that
\begin{equation}
    \limsup\nolimits_{n\to\infty}
    \E[\|\Theta_n-\vartheta\|^2]
    =
    0.
\end{equation}}
\end{aproof}

\subsubsection*{Acknowledgements}
This work has been supported by the Ministry of Culture and Science NRW as part
of the Lamarr Fellow Network. In addition, this work has been partially funded by the European Union (ERC, MONTECARLO, 101045811). The views and the opinions expressed in this work are however those of the authors only and do not necessarily reflect those of the European Union or the European Research Council (ERC). Neither the European Union nor the granting authority can be held responsible for them. Moreover, we gratefully acknowledge the Cluster of Excellence EXC 2044/2-390685587, Mathematics Münster: Dynamics-Geometry-Structure funded by the Deutsche Forschungsgemeinschaft (DFG, German Research Foundation).
Most of the specific formulations in the proofs of this work have been created using \cite{kuckuck2025some}.

\subsubsection*{Use of large language models}
GPT-5.6 Pro has supported in extending and revising the arguments in some of the proofs of this work and has supported us in creating the literature review in \cref{lite}.
The authors take full responsibility for each of the sentences/statements made in the paper.

{\small
\bibliography{ref}
\bibliographystyle{acm}}

\end{document}